\documentclass[a4paper,10pt]{amsart}

\usepackage[utf8]{inputenc}
\usepackage[T1]{fontenc}        

\usepackage{hyperref}
\usepackage{amsmath}
\usepackage{amsfonts}
\usepackage{amssymb}
\usepackage{amsthm}
\usepackage{graphicx}
\usepackage{geometry}
\usepackage{float}

\usepackage{tikz}
\usepackage{pgfplots}
\usetikzlibrary{pgfplots.groupplots}
\usetikzlibrary{matrix}     % Pour faire des diagrammes

\newtheorem{lem}{Lemma}
\newtheorem{thrm}{Theorem}
\newtheorem{prop}{Proposition}

\newtheorem{rem}{Remark}
\newtheorem{example}{Example}
\newtheorem{cor}{Corollary}

\numberwithin{equation}{section}
\newcommand{\ds}{\displaystyle}

\newcommand{\pa}{\partial}
\newcommand{\R}{\mathbb R}
\newcommand{\T}{\mathbb{T}}
\newcommand{\I}{\mathcal{I}}
\newcommand{\J}{\mathcal{J}}
\newcommand{\Jhalf}{\mathcal{J}^*}
\newcommand{\Z}{\mathbb{Z}}
\newcommand{\F}{\mathcal{F}}
\newcommand{\G}{\mathcal{G}}
\newcommand{\M}{\mathcal{M}}

\newcommand{\HHH}{\mathbb{H}}

\newcommand{\V}{\mathcal{V}}
\newcommand{\X}{\mathcal{X}}
\newcommand{\U}{\mathbf{U}}

\newcommand{\cH}{\mathcal{H}}
\newcommand{\lla}{\left\langle}
\newcommand{\rra}{\right\rangle}
\newcommand{\dD}{\mathrm{d}}
\newcommand{\ve}{\varepsilon}
\newcommand{\eg}{\emph{e.g. }}

\newcommand{\RR}{\mathbb{R}}
\newcommand{\NN}{\mathbb{N}}
\newcommand{\MM}{\mathcal{M}}

\newcommand{\TT}{\mathbb{T}}
\newcommand{\CC}{\mathbb{C}}
\newcommand{\Q}{\mathcal{Q}}

\newcommand\Tom[1]{\textcolor{black}{#1}}
\newcommand\Mar[1]{\textcolor{black}{#1}}
\newcommand\Max[1]{\textcolor{black}{#1}}
\newcommand\mysubsection[1]{\vspace{1em}\textbf{#1}}

\title[Numerical hypocoercivity and diffusion limit]{Hypocoercivity and diffusion limit of a finite volume scheme for linear kinetic equations}

\author{Marianne Bessemoulin-Chatard}%
\address{\flushleft
  \underline{Marianne Bessemoulin-Chatard}\\[.3em]
  Université de Nantes, Laboratoire de Mathématiques Jean Leray, UMR CNRS 6629\\[.2em]
  F-44000 Nantes, France\\[.2em]
  Email: \texttt{marianne.bessemoulin@univ-nantes.fr}}%
\thanks{MBC was partially funded by the Centre Henri Lebesgue (ANR-11-LABX-0020-01) and ANR Project MoHyCon (ANR-17-CE40-0027-01).}%

\author{Maxime Herda}
\address{\flushleft
  \underline{Maxime Herda}\\[.3em]
  Inria, Univ. Lille, CNRS, UMR 8524 – Laboratoire Paul Painlevé\\[.2em]
  F-59000 Lille, France\\[.2em]
  Email: \texttt{maxime.herda@inria.fr}}
\thanks{MH  was partially funded by ANR Project MoHyCon (ANR-17-CE40-0027-01).}

\author{Thomas Rey}
\address{\flushleft
  \underline{Thomas Rey}\\[.3em]
  Univ. Lille, CNRS, UMR 8524, Inria – Laboratoire Paul Painlevé\\[.2em]
  F-59000 Lille, France\\[.2em]
  Email: \texttt{thomas.rey@univ-lille.fr}}
\thanks{TR were partially funded by Labex CEMPI (ANR-11-LABX-0007-01) and ANR Project MoHyCon (ANR-17-CE40-0027-01). }

\begin{document}
  
  \maketitle
  
  \begin{abstract}
    In this article, we are interested in the asymptotic analysis of a finite volume scheme for one dimensional linear kinetic equations, with either Fokker-Planck or linearized BGK collision operator. Thanks to appropriate uniform estimates, we establish that the proposed scheme is Asymptotic-Preserving in the diffusive limit. Moreover, we adapt to the discrete framework the hypocoercivity method proposed by [J. Dolbeault, C. Mouhot and C. Schmeiser, \textit{Trans. Amer. Math. Soc.}, 367, 6 (2015)] to prove the exponential return to equilibrium of the approximate solution. We obtain decay \Max{rates that are bounded uniformly} in the diffusive limit.
    Finally, we present an efficient implementation of the proposed numerical schemes, and perform numerous numerical simulations assessing their accuracy and efficiency in capturing the correct asymptotic behaviors of the models.\\[1em]
    \textsc{Keywords:} Kinetic equations, finite volume methods, hypocoercivity, diffusion limit, asymptotic-preserving schemes.\\[.5em]
    \textsc{2010 Mathematics Subject Classification:} 82B40, %Kinetic theory of gases
    65M08, %Finite volume methods
    65M12. %Stability and convergence of num. methods
  \end{abstract}
  
  \tableofcontents
  
  \section{Introduction}
  
  \textbf{The linear kinetic equation.} Many engineering or biological problems involve fluid-like systems in transitional regimes: micro-electro-mechanical systems, space shuttle reentry, powder and grains in storage silos, bacteria colonizing a medium, \emph{etc}. As a consequence of nonequilibrium behaviors, the macroscopic descriptions can break down, and a kinetic model may be needed to depict accurately the system. In this contribution, we are interested in the numerical approximation of a prototypical model for such systems. It describes the evolution of a particle distribution function $f^\ve = f^\ve(t,x,v)$, for time $t \geq 0$, at position $x \in \mathbb{T}$ the one dimensional torus, and velocity $v \in \RR$, solution to the initial value problem
  \begin{equation} 
    \label{eq:Collision}
    \left\{ \begin{aligned}
      & \ve \frac{\partial f^\ve}{\partial t} + v \frac{\partial f^\ve}{\partial x}\,=\, \frac1\ve\,\Q(f^\ve), 
      \\ & \, \\
      & f^\ve(0, x, v) = f_{0}(x,v) \geq 0\,.
    \end{aligned} \right.
  \end{equation}
  The \emph{collision} operator $\Q$ describes a microscopic ``collision process'' acting only on the velocity variable $v$ and preserving the zeroth order moment of the distribution $f^\ve$, namely the mass or total number of particles. The small scaling parameter $\ve > 0$ is the ratio between the mean free path of particles  and the length scale of observation. In the context of rarefied gas dynamics it corresponds to the dimensionless \emph{Knudsen} number. The presence of the factor $\ve$ in front of the time derivative accounts for the fact that the system is observed on long time scales. The \Mar{system} is said to be in the kinetic regime if $\ve \sim 1$ and in the diffusive regime if $\ve \ll 1$.
  
  Throughout this paper, the collision operator will be either of linear Fokker-Planck type:
  \begin{equation}
    \label{def:Qfp}
    \Q_{FP}(f)(v) = \partial_v\left (\partial_v f + vf \right )
  \end{equation}
  or of linearized BGK \cite{Bhatnagar:1954} / relaxation type:
  \begin{equation}
    \label{def:Qbgk}
    \Q_{BGK}(f)(v) = \rho M(v) - m_0 f(v),\qquad\rho = \int_\RR f(v) \, \dD v
  \end{equation}
  where $M$ is a given nonnegative function of $v$, the so-called \emph{Maxwellian} whose moments are denoted by
  \[
  m_k\ =\ \int_\RR |v|^k M(v)\,\dD v.
  \]
  In particular, $m_0$ is the mass of $M$.
  
  \begin{rem}
    The linear BGK model appears quite naturally as a linearization of the complete, quadratic, Boltzmann equation \cite{CIP:94} near its global Maxwellian equilibrium $M$. Indeed, one can write 
    \begin{align*}
      \Q_{BGK}(f)(v) & = \Mar{M(v)}\int_\RR f(v_*) \, \dD v_* -f(v) {\int_\RR M(v_*) \, \dD v_*} \\
      & = \int_\RR \left [ f_* \, M - f \, M_* \right ] \dD v_* ,
    \end{align*}
    where  the usual shorthand notations $\psi := \psi(v)$ and $\psi_* := \psi(v_*)$ were used. 
  \end{rem}

  In the case of the BGK operator, every result in this article \Mar{applies} to Maxwellians $M$ that are even functions of $v$ and which have finite moments $m_k<+\infty$ up to $k = 4$ \Max{(see Section~\ref{s:nongauss} for an illustration)}.  However, in order to unify the analysis of both collision operators, we will focus in the rest of the presentation on the classical Boltzmann's Maxwellian, namely the centered reduced Gaussian function
  \[
  M(v) := \frac 1{\sqrt{2 \pi}} \exp\left (-\frac{|v|^2}{2}\right ).
  \]
  In particular, $m_0 = 1$ and $\Q_{BGK}(f) = \rho M - f$.

  \begin{rem}
    Observe that such a Maxwellian distribution verifies the elementary identity
    \begin{equation}
      \label{eq:MaxwellianIdentity}
      M'(v) = -v M(v).
    \end{equation}
    Hence, one can rewrite the Fokker-Planck operator in the following gradient form
    \begin{equation}
      \label{def:Qfp_grad}
      \Q_{FP}(f) = \partial_v \left (M \, \partial_v \left (\frac{f}{M}\right )\right ).		  
    \end{equation}
  \end{rem}
  
  The motivations in this paper are twofold and \Mar{concern} the design and analysis of a numerical scheme which is accurate in two asymptotics of \eqref{eq:Collision}. The first \Tom{asymptotics} is the diffusion limit $\ve\to0$ and the second is the long-time \Tom{asymptotics} $t\to\infty$.
  
  \mysubsection{Diffusion limit.}
  A natural mathematical problem concerning kinetic equations such as \eqref{eq:Collision} is the study of the diffusion / parabolic limit $\ve\to0$. Thanks to the  choice of time scale, it is possible to capture nontrivial macroscopic dynamics at the limit. In order to identify the limit, one can multiply the equation \eqref{eq:Collision} by $(1,v)$ and integrating in the velocity space to get
  \begin{equation}
    \label{eq:momentsCoupled}
    \left \{ \begin{aligned}
      & \partial_t \rho^\ve + \partial_x j^\ve = 0, \\
      & \ve^2 \partial_t j^\ve +  m_2\partial_x \rho^\ve +\partial_x S^\ve  = - j^\ve,
    \end{aligned} \right.
  \end{equation}
  where the moments are defined by
  \begin{equation*}
    \rho^\ve := \int_\RR f^\ve \, \dD v, \quad j^\ve := \frac 1\ve \int_\RR f^\ve \, v \, \dD v, \quad S^\ve := \int_\RR (v^2 - m_2) f^\ve \, \dD v.
  \end{equation*}
  Then, by formally taking limits $\ve \to 0$ in \eqref{eq:Collision} and \eqref{eq:momentsCoupled}, one obtains $f^\ve\to f = \rho\,M$ \Max{which implies}
  $
  S^\ve \to 0$ \Max{and the second equation of \eqref{eq:momentsCoupled} yields}
 $
 j^\ve \to - m_2 \partial_x \rho
  $.  Hence, the limit density $\rho$ is solution to the linear heat equation 
  \begin{equation}
    \label{eq:heatEQ}
    \frac{\partial \rho}{\partial t} - \partial_x (m_2 \partial_x \rho )=0.
  \end{equation}  
  
  In the literature, the first results concerning the approximation of kinetic models by diffusion equations have been proposed in \cite{bensoussan1979boundary,Keller1974}. Justifications of asymptotic expansion of the solutions in power of $\ve$ can be found in \cite{Bardos1984} for the neutron transport, in \cite{Degond1987} for the Fokker-Planck equation, in \cite{Bardos1988} for the radiative transfer equation, and in \cite{poupaud1991diffusion} for linear semiconductor Boltzmann equation. In \cite{Lions1997}, the authors studied the diffusive limit of generalized two-velocity kinetic models. A large class of linear collision operators was considered in \cite{Degond2000a}, and approximation of kinetic equations by diffusion is justified by homogenization in \cite{Goudon2001}. We also refer to \cite[Chapter~2]{bouchut_2000_kinetic} for an overview of classical hydrodynamic and parabolic limits of kinetic equations. Extensions to more general models with fractional diffusion limits were more recently obtained in \cite{mellet2011fractional}.
  
\begin{figure}
    \label{fig:APdiag}
    \tikzstyle{block} = [rectangle, draw, text width=2em, text centered, rounded corners, minimum height=2em]
    \tikzstyle{line} = [draw, -latex]
    \begin{center}
      \begin{tikzpicture}[node distance=6em]
	
	\node [block, text width=13em, minimum height=10em] (fgraph) {};
	
	\node [left of=fgraph, xshift=3em, yshift=2.5em] (fepsdel) {\large$f_h^{\ve}$};
	\node [right of=fepsdel] (feps) {\large$f^{\ve}$};
	\node [below of=fepsdel, node distance=5em](fdel) {\large$f_h$};
	\node [right of=fdel]  (f) {\large$f$};
	eq:
	\path        (fepsdel) -- node[above]{\small$h\to0$}(feps);
	\path [line] (fepsdel) -- node[left, yshift=.3em]{\small$\ve\to0$}(fdel);
	\path        (feps) -- node[right, yshift=.3em]{\small$\ve\to0$}(f);
	\path [line] (fdel) -- node[below]{\small$h\to0$}(f);
	
	\draw[dashed,->] (fepsdel) -- (feps);
	\draw[dashed,->] (feps) -- (f);
	
      \end{tikzpicture}
    \end{center}
    \caption{The AP diagram ($h$ denotes the size of the discretization)}
  \end{figure}
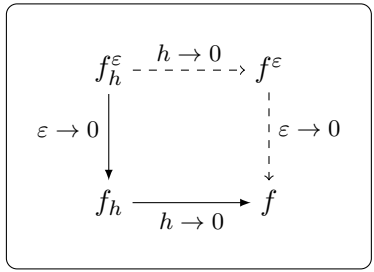

   There exists a dense literature describing numerical schemes which enjoy the property of being stable in the diffusion limit, and converge to a numerical discretization of  the macroscopic model. These schemes fall down in the so-called Asymptotic-Preserving (AP) framework and are compatible with the kinetic and asymptotic regime of the equation. The principle of AP schemes can be roughly summarized as the commutative diagram presented in Fig. \ref{fig:APdiag}.
  
  In the framework of kinetic equations, AP schemes first appeared 20 years ago, \Mar{on the one hand for kinetic semiconductors models in two papers due to Klar \cite{klar1998asymptotic} for the linear BGK equation and to Schmeiser and Zwirchmayr \cite{Schmeiser1999} for the \Tom{linear} Boltzmann equation, and on the other hand for relaxation models (\eg two velocities Goldstein-Taylor) in the article of Jin, Pareschi, and Toscani \cite{jin2000uniformly}.}
  Two recent review papers \cite{jin:2010portoercole, DimarcoPareschi15} contain an almost up to date bibliography on numerical methods for collisional kinetic equations of type \eqref{eq:Collision} and AP schemes.
  For even more recent works on AP schemes, one can consult these two papers \cite{dimarco2018asymptotic, CrestettoCrouseillesLemou2017} on \Mar{particle} methods and the references therein.
  
    Nevertheless, most of these works contain only formal proofs of the AP properties, the very important solid lines of Fig. \ref{fig:APdiag}.
  Here are some instances of AP schemes with proofs: the micro-macro approach of Lemou and Mieussens \cite{LemMieu:2008}, for either discrete velocity models (Goldstein-Taylor-like) or linear scattering  model (continuous in velocity), \Mar{contain} some partial proofs. These proofs were later completed by Liu and Mieussens in \cite{liu2010analysis} with a very detailed stability analysis. Note that the numerical scheme for dealing with the kinetic equation is semi-implicit in time, but the limiting  $\ve \to 0$ scheme is explicit, yielding a very restrictive parabolic CFL condition on the time step.
  For the anomalous diffusion limit, let us also mention the work of Crouseilles, Hivert, and Lemou \cite{crouseilles2016numerical} where an AP scheme is built, with a proof of convergence (semi-discrete in time), without rate of convergence towards the equilibrium.

   \mysubsection{Hypocoercivity and large time behavior.} 
   A second classical problem in the asymptotic analysis of kinetic equations is the large time behavior of solutions. On the class of models under consideration, there are several ways to show that there are constants $\kappa>0$ and $C\geq1$ such that the solution satisfies
   \begin{equation}
     \| f^\ve(t) - \mu_f M \|_\mathcal{X} \leq C \| f_0 - \mu_f M \|_\mathcal{X} e^{-\kappa t},
     \label{eq:decay}
   \end{equation}
   where \[\mu_f=\iint_{\TT\times\RR} f_0 \,\dD x\,\dD v\,,\]
   and $\mathcal{X}$ is some appropriate functional space. \Max{In the diffusive scaling under consideration}, \Mar{the constants $C$ and $\kappa$ can be chosen independently of $\varepsilon$ as $\varepsilon$ tends to zero.}
   In the case of the kinetic Fokker-Planck operator, the model is simple enough so that it is actually possible to compute explicitly the fundamental solution and derive such an estimate. 
   
   However, in recent years, more systematic and robust methods have been developed in order to show \eqref{eq:decay}. They are called \emph{hypocoercivity} methods and can usually be seen as entropy methods \cite{jungel_2016_entropy, schmeiser2018entropy}, in the sense that the goal is to design a functional $H(f^\ve)$, called the (modified) entropy, which is dissipated along the solutions. From there a control of the entropy by its dissipation provides decay of the former with explicit rates. Hypocoercivity methods deal with the fact that for kinetic \Mar{equations} such as \eqref{eq:Collision}, dissipation / relaxation apparently occurs only in the velocity variable. Nonetheless, thanks to the mixing properties of the transport operator $v\partial_x$ in the phase space $x-v$, one can recover dissipation along $x$ and $v$. The denomination hypocoercivity is to be related with the fact that if \eqref{eq:decay} held with $C=1$ then the operator $v\partial_x - \Q(\cdot)$ would be coercive in the usual sense.
   
  A comprehensive introduction on hypocoercivity methods may be found in the lectures of Hérau \cite{herau_2018_introduction}. The first hypocoercivity method introduced by Hérau and Nier for kinetic Fokker-Planck equations \cite{herau_2004_isotropic} is very similar to the method of commutators for proving hypoellipticity going back to the seminal work of Hörmander \cite{hormander1967}. This method has been successfully applied to non-hypoelliptic  equations (without regularization effects), such as linear Boltzmann or BGK equations, by Hérau \cite{herau_2006_hypocoercivity}, and Mouhot and Neumann \cite{mouhot_2006_quantitative}. A general abstract framework for this method of commutators has then been given by Villani in \cite{Villani:2009Hypo}. The drawback of these hypocoercivity methods is that they often require regularity on the initial datum and involve estimates in equivalent weighted $H^1$ norms leading to sometimes tedious computations. \Max{In \cite{herau_2006_hypocoercivity}, the technique is adapted to start from only weighted $L^2$ initial datum, but the structure of the functional is similar.}
  
      In a recent work \cite{Dolbeault2015}, Dolbeault, Mouhot and Schmeiser established a new hypocoercivity method in a general abstract framework. It is based on the definition of a modified entropy taking the form of an equivalent weighted $L^2$ norm \Max{involving macroscopic quantities (\emph{i.e.} integrated in velocity)}. In this work we focus on this particular hypocoercivity method\Max{, and we show that it is well-suited for an adaptation in the discrete framework.}

  Concerning the numerical methods preserving large-time behaviors of solutions, several efficient numerical schemes have been developed for homogeneous kinetic equations. A full discrete finite difference scheme for the homogeneous Fokker-Planck equation was built in the pioneering work of Chang and Cooper \cite{ChangCooper:1970}. More recently, several schemes preserving the exponential trend to equilibrium have been proposed for nonlinear degenerate parabolic equations (see for example \cite{Bessemoulin-Chatard2012,Burger2010,Chainais-Hillairet2007,Gosse2006}). In \cite{Pareschi2017}, the question of describing correctly the equilibrium state of several nonlinear diffusion and kinetic equations is addressed together with that of the order of the schemes. Let us also mention \cite{Filbet2017,chainais_2018_large} where the case of non-homogeneous Dirichlet boundary conditions are dealt with. 
  
  In the case of inhomogeneous kinetic equations, only few results are available, as the understanding of hypocoercive structures is quite recent. The case of the very simple Kolmogorov equation has been investigated in \cite{poretta_2017_numerical,Foster2017,Georgoulis}. In \cite{Foster2017}, a time-splitting technique based on self-similarity properties is used for solutions that decay polynomially in time. In \cite{poretta_2017_numerical,Georgoulis}, the $H^1$ method \emph{à la} Hérau and Villani has been adapted for respectively a finite difference and a finite element schemes. Lastly, let us finish with the recent paper of Dujardin, Hérau and Lafitte \cite{Dujardin} where the same $H^1$ method is adapted to a finite difference scheme for the kinetic Fokker-Planck equation \Mar{\eqref{eq:Collision}--\eqref{def:Qfp}}.
  
  \mysubsection{Main results.} In this contribution we design implicit in time finite volume schemes for \eqref{eq:Collision} with collision operators \eqref{def:Qfp} and \eqref{def:Qbgk}. We perform their numerical asymptotic analysis and are able to show the following properties.
  \begin{itemize} 
    \item Our schemes are AP (Theorem \ref{th:AP}) in the diffusion limit $\ve\to0$. Solutions converge to solutions of a finite difference scheme solving the heat equation \eqref{eq:heatEQ}. Provided that the scheme is implemented correctly (see Section~\ref{sec:implementation}) one can take $\ve=0$ in the scheme. In particular the linear system that has to be solved at each time iteration has coefficients and condition number that are uniformly bounded with respect to $\ve$.
    \item The solutions of our schemes preserve the decay of a discrete analogue of the modified entropy of Dolbeault, Mouhot and Schmeiser \cite{Dolbeault2015}. Consequently, we are able to prove the discrete equivalent of the exponential decay estimate \eqref{eq:decay} (Theorem~\ref{th:hypocoercivity_discrete}). The decay rates are bounded uniformly in the diffusion limit.
  \end{itemize}
  Some comments are in order. On the first hand, let us mention that an essential identity which allows to prove both properties in the Fokker-Planck case is a discrete equivalent of \eqref{eq:MaxwellianIdentity}. By choosing the framework of finite volume schemes we can introduce a discrete Maxwellian at cell centers and interfaces of the mesh in velocity. Enforcing the discrete equivalent of \eqref{eq:MaxwellianIdentity} just relates these two quantities and leaves a large choice in the actual definition of the discrete Maxwellian.

  On the second hand, by choosing to adapt the $L^2$ hypocoercivity method rather than the $H^1$ method, we avoid the burden of estimating derivatives of the solution and dealing with commutators coming from the choice of discrete derivative operators. Thus, the numerical analysis turns out to be almost identical to what is done in the continuous case. For this reason, we are confident that this method could successfully be applied for the asymptotic analysis of other numerical schemes for hypocoercive equations.

  \mysubsection{Plan of the paper.} The outline of this article is as follows. In Section \ref{sec:Continuous}, we briefly present the results in the continuous setting. We start by establishing some uniform estimates, which are the key point to study both the diffusive limit and the exponential return to equilibrium. Then we state the convergence of the kinetic equation towards the heat equation as $\ve$ tends to 0. To conclude this section, we prove the $L^2$-hypocoercivity of both Fokker-Planck and linearized BGK operators, by introducing a modified entropy functional which is a slight simplification of the original version proposed by Dolbeault, Mouhot and Schmeiser \cite{Dolbeault2015}. In Section \ref{sec:FullDiscrete}, we adapt these results to the discrete framework. The considered numerical schemes are implicit in time, and finite volume in space and velocity. We discretize carefully the Fokker-Planck fluxes, in order to obtain a discrete version of the $L^2$ entropy estimate, which is the starting point to establish the needed discrete uniform estimates. Then our two main results are stated in Theorem \ref{th:AP} (asymptotic-preserving property in the diffusive limit) and in Theorem \ref{th:hypocoercivity_discrete} (numerical hypocoercivity). Section \ref{sec:implementation} is then devoted to the practical implementation of our method. We present a discretization equivalent to the proposed schemes, based on a perturbative micro-macro formulation, and we write the obtained linear systems in explicit matrix forms for sake of completeness. Finally in Section \ref{sec:numres}, we present several numerical experiments which demonstrate the efficiency of our schemes. After investigating numerically the asymptotic-preserving property, we finally study the trend to equilibrium for both Fokker-Planck and linearized BGK operators.

  \mysubsection{Acknowledgments.} The authors would like to thank G. Dujardin, F. Hérau and P. Lafitte for fruitful discussions during the conception of this article.
  
\section{The continuous setting} 
\label{sec:Continuous}

For convenience, let us introduce the measures
$
\dD M  := M(v)\dD v\,,
$
\[
\dD \gamma  = \gamma(v)\dD v := \tfrac{\dD v}{M(v)}
\]
and $L^2(\dD M)$ (\textit{resp.} $L^2(\dD\gamma)$) the spaces of square integrable functions against the measure $\dD M$ (\textit{resp.} $\dD\gamma$). A straightforward consequence of the Gaussian Poincaré inequality is that for any $f\in L^2(\dD\gamma)$, one has
 \begin{equation}
 \|f-\rho M\|_{L^2(\dD\gamma)}\,\leq\,\|\partial_v\left(\tfrac{f}{M}\right)\|_{L^2(\dD M)}\,,
 \label{eq:poincare}
 \end{equation}
 where $\rho = \int f\dD v$.

\subsection{Uniform estimates}\label{sec:unif}

With an initial data in $L^2(\dD x\dD \gamma)$, there is a unique solution to \eqref{eq:Collision}, with either Fokker-Planck or BGK collision operators (see \eg \cite{AllaireBlancDespresGolse:lecturenotes}). This solution conserves mass and nonnegativity. Additionally one has the following basic ``entropy'' estimate.

\begin{lem}
 Let $f^\ve$ be the solution of \eqref{eq:Collision}  with initial data $f_0\in L^2(\dD x\dD \gamma)$. Then for every $t\geq 0$\,,

 \begin{equation}
 \frac{1}{2} \frac{\dD}{\dD t}\|f^\ve(t)\|_{L^2(\dD x\dD\gamma)}^2\,+\,\frac{1}{\ve^2}\|f^\ve(t)-\rho^\ve(t) M\|^2_{L^2(\dD x\dD\gamma)}\ \leq\ 0\,.
  \label{eq:entropy1}
 \end{equation}
 In particular one has
 \begin{equation}
  \max\left(\|f^\ve\|_{L^\infty_tL^2(\dD x\dD\gamma)}^2, \frac{2}{\ve^2}\|f^\ve - \rho^\ve M\|_{L^2_tL^2(\dD x\dD\gamma)}^2\right)\,\leq\,\|f_0\|_{L^2(\dD x\dD\gamma)}^2\,.
  \label{eq:L2gamma}
 \end{equation}

\label{lem:basic_estimates}
\end{lem}
\begin{proof} Estimate \eqref{eq:entropy1} is obtained by multiplying \eqref{eq:Collision} with $f^\ve/M$ and integrating in $x$ and $v$ yielding
 \[
   \frac12\frac{\dD}{\dD t}\|f^\ve\|^2_{L^2(\dD x\dD\gamma)} + \frac{1}{2\ve}\iint \partial_x\left(v\left(\frac{f^\ve}{M}\right)^2\right)\dD x\dD v\ = \ \frac{1}{\ve^2}\iint \Q(f^\ve)\frac{f^\ve}{M}\dD x\,\dD v\,.
 \]
 The second term of the left hand side cancels. \newline \noindent In the case of the BGK operator, since $\int \Q_{BGK}(f^\ve)\dD v= 0$, one has
 \[
 \iint \Q_{BGK}(f^\ve)\frac{f^\ve}{M}\dD x\,\dD v\ =\ \iint \Q_{BGK}(f^\ve)\left(\frac{f^\ve}{M}-\rho^\ve\right)\dD x\,\dD v\ =\ - \|f^\ve - \rho^\ve M\|_{L^2(\dD x\dD\gamma)}^2\,.
 \]
 In the case of the Fokker-Planck operator, an integration by parts in $v$ and the Poincaré inequality \eqref{eq:poincare} yields
  \[
 \iint \Q_{FP}(f^\ve)\frac{f^\ve}{M}\dD x\,\dD v\ =\ - \|\partial_v\left(\frac{f^\ve}{M}\right)\|_{L^2(\dD x\dD M)}^2\ \leq\ -\|f^\ve-\rho^\ve M\|_{L^2(\dD x\dD\gamma)}^2.
 \]
\end{proof}

From the entropy estimate, one can gather many estimates on the moments.

\begin{lem}[Moments estimates]
 Under the hypotheses of Lemma~\ref{lem:basic_estimates}, the moments $\rho^\ve$, $j^\ve$ and $S^\ve$ satisfy equations \eqref{eq:momentsCoupled} and  the following estimates:
 \begin{equation}
  \|\rho^\ve(t)\|_{L^2_x}\ \leq\ \|f^\ve(t)\|_{L^2(\dD x\dD\gamma)}\,,
  \label{eq:estmomentsLinfrho}
 \end{equation}
 \begin{equation}
 \ve\,\|j^\ve(t)\|_{L^2_x}\ \leq\ m_2^{1/2}\,\|f^\ve(t)\|_{L^2(\dD x\dD\gamma)}\,
  \label{eq:estmomentsLinfj}
 \end{equation}
 \begin{equation}
 \ve\,\|j^\ve(t)\|_{L^2_x}\ \leq\ m_2^{1/2}\,\|f^\ve(t)-\rho^\ve(t)\,M\|_{L^2(\dD x\dD\gamma)}\,,
  \label{eq:estmomentsL2j}
 \end{equation}
 \begin{equation}
 \|S^\ve(t)\|_{L^2_x}\ \leq\ (m_4 - m_2^2)^{1/2}\,\|f^\ve(t)-\rho^\ve(t)\,M\|_{L^2(\dD x\dD\gamma)}\,.
  \label{eq:estmomentsL2S}
 \end{equation}
 In particular, there is a constant $C>0$ depending only on $m_2$ and $m_4$ such that
 \begin{equation}
\max\left(\|\rho^\ve\|_{L^\infty_tL^2_x},\,\ve\,\|j^\ve\|_{L^\infty_tL^2_x},\,\|j^\ve\|_{L^2_tL^2_x},\, \frac1\ve\,\|S^\ve\|_{L^2_tL^2_x}\right)\ \leq\ C\,\|f_0\|_{L^2(\dD x\dD\gamma)}\,.
\label{eq:estimmoments}
 \end{equation}
\label{lem:moments}
\end{lem}
\begin{proof}
 In order to estimate the moments, let us observe that one can rewrite the macroscopic quantities as
 \begin{align*}
     & \rho^\ve=\int M^{1/2}(f^\ve M^{-1/2}) \,\dD v, \\
     & \ve\,j^\ve = \int (v M^{1/2})(f^\ve M^{-1/2})\,\dD v\ =\ \int (v M^{1/2})(f^\ve-\rho^\ve M) M^{-1/2}\,\dD v, \\
     & S^\ve = \int (v^2 - m_ 2)M^{1/2}(f^\ve-\rho^\ve M) M^{-1/2}\,\dD v,
 \end{align*}
 by symmetry of $M$ for $j^\ve$, and by definition of $m_2$ for $S^\ve$. 
 Then one can use Cauchy-Schwarz inequalities to get the desired estimates. Estimate \eqref{eq:estimmoments} is then a consequence of \eqref{eq:L2gamma}. Finally, integrate \eqref{eq:Collision} in $v$ against $1$ and $v$ to get \eqref{eq:momentsCoupled}. 
\end{proof}

\subsection{Diffusive limit}
\begin{prop}
Let $(f^\ve)_{\ve\geq0}$ be a sequence of solutions of \eqref{eq:Collision} (with either Fokker-Planck or BGK collision operator) with initial data $f_0\in L^2(\dD x\dD \gamma)$. Then when $\ve\to0$, one has
\[
 f^\varepsilon(t,x,v)\rightarrow f(t,x,v)\ =\ \rho(t,x)\,M(v)\;\text{weakly in }L^2_{t,\text{loc}}L^2(\dD x\dD\gamma) 
\]
and $\rho$ solves the heat equation
\[
 \partial_t\rho - m_2\,\partial_{xx}^2\rho\ =\ 0\,,\quad\rho(0,x)\,=\,\int f_0(x,v)\,\dD v.
\]
\end{prop}
\begin{proof}
 Using that $f^\ve - \rho^\ve M$ is uniformly $O(\ve)$ in $L^2(0,\infty;L^2(\dD x\dD\gamma))$ and that $\rho^\ve$ is uniformly bounded in $L^\infty_tL^2_x$, we get the existence of the limit, up to extraction of a subsequence. Then, by combining both equations of \eqref{eq:momentsCoupled} one gets
 \[
  \partial_t\rho^\ve - m_2\partial_{xx}^2\rho^\ve\ =\ \ve^2\partial_tj^\ve + \partial_xS^\ve
 \]
 in the sense of distributions. Taking limits thanks to the estimates \eqref{eq:estimmoments} one recovers the heat equation. By uniqueness of the limit the whole sequence converges.
 \end{proof}

\subsection{Hypocoercivity}\label{sec:hypocoercivity}

Let $\cH$ be the subspace of the Hilbert space $L^2(\dD x\dD\gamma)$ composed of the zero mean functions with canonical scalar product. The space $\cH$ \Mar{coincides} with the orthogonal space of the nullspace of the unbounded operator $L_\ve = \tfrac1\varepsilon v\partial_x - \tfrac{1}{\ve^2}\Q(\cdot)$ in both the Fokker-Planck and BGK cases. One can show thanks to the method developed in \cite{Dolbeault2015} that the latter operator is hypocoercive (see \cite{Villani:2009Hypo}) in $\cH$ uniformly for small $\ve$. More precisely, one has the following result.
\begin{prop}
 There are constants $C\geq 1$ and $\kappa>0$ such that for all $\ve\in(0,1)$ and all initial data $f_0\in L^2(\dD x\dD \gamma)$,  the solution $f^\ve$ of \eqref{eq:Collision} satisfies
 \begin{equation*}
 \|f^\ve(t)-\mu_f M\|_{L^2(\dD x\dD \gamma)}\ \leq\ C\,\|f_0-\mu_f M\|_{L^2(\dD x\dD \gamma)}\,e^{-\kappa\,t}.
 \end{equation*}
 \label{prop:hypocoercivity}
\end{prop}
 The rest of the section is devoted to the proof of Proposition~\ref{prop:hypocoercivity}. We need some intermediate results stated in the following lemmas. Up to changing $f_0$ by $f_0-\mu_f M$ and correspondingly $f^\ve$ and $\rho^\ve$ by respectively $f^\ve -\mu_f M$ and $\rho^\ve -\mu_f$, we can assume from now on in this section that
 \[
  \mu_f\ =\ \iint f_0 \dD x\dD v\ =\ 0\,.
 \]
It follows that $f^\ve$ and $\rho^\ve$ are mean-free.

 In order to prove Proposition~\ref{prop:hypocoercivity}, we introduce the ``modified entropy functional'' \textit{à la }Dolbeault, Mouhot and Schmeiser. It is a slight simplification of the original version of \cite{Dolbeault2015} using the fact that we work on a bounded space domain. It reads
\begin{equation}
 H(f^\ve(t)) = \frac12\|f^\ve(t)\|_{L^2(\dD x\dD \gamma)}^2\,+\,\eta\,\ve^2\lla j^\ve, \partial_x\phi^\ve\rra_{L^2_x}\,,
 \label{eq:modifiedEntropy}
\end{equation}
where $\phi^\ve(t,x)$ is the solution of the Poisson equation
\[
 -\partial_{xx}^2\phi^\ve\ =\ \rho^\ve\,,\qquad \int_\mathbb{T}\phi^\ve\dD x\,=\,0
\]
and $\eta>0$ is a small parameter to be chosen. 
\begin{lem}\label{lem:estimphi}
 The function $\phi^\ve$ satisfies for all $t\geq0$
 \begin{eqnarray*}
  \|\partial_x\phi^\ve(t)\|_{L^2_x}&\leq&C_P\,\|\rho^\ve(t)\|_{L^2_x},\\
  \|\partial_t\partial_x\phi^\ve(t)\|_{L^2_x}&\leq&\|j^\ve(t)\|_{L^2_x},
 \end{eqnarray*}
 where $C_P = 1/2\pi$ is the Poincaré constant of $\mathbb{T}$.
\end{lem}
\begin{proof}
 The first estimate is obtained by testing the Poisson equation against $\phi^\ve$ and applying the Poincaré inequality
 \[
  \|\partial_x\phi^\ve\|_{L^2_x}^2 \,=\, \lla-\partial_{xx}^2\phi^\ve, \phi^\ve\rra_{L^2_x}\,\leq\,\|\rho^\ve\|_{L^2_x}\,\|\phi^\ve\|_{L^2_x}\,\leq\,C_P\|\rho^\ve\|_{L^2_x}\,\|\partial_x\phi^\ve\|_{L^2_x}\,.
 \]
 For the second estimate, first take the time derivative of the Poisson equation and use the continuity equation $\partial_t\rho^\ve + \partial_xj^\ve = 0$ to get
 $
  -\partial_t\partial_{xx}^2\phi^\ve\ =\ -\partial_xj^\ve
 $. 
Then multiply by $\partial_t\phi^\ve$ and integrate to get 
\[
 \|\partial_t\partial_x\phi^\ve\|_{L^2_x}^2\, =\, \lla-\partial_xj^\ve, \partial_t\phi^\ve\rra = \lla j^\ve, \partial_t\partial_x\phi^\ve\rra_{L^2_x}\,\leq\,\|j^\ve\|_{L^2_x}\,\|\partial_t\partial_x\phi^\ve\|_{L^2_x}.
\]
\end{proof}

For small enough $\eta>0$, the square root of the modified entropy actually defines an equivalent norm on $\cH$.

\begin{lem}[Equivalent norm]
 There is $\eta_1\Mar{>0}$ such that for all $\eta\in(0,\eta_1)$, there are positive constants $0<c_\eta<C_\eta$ such that for all $\ve\in(0,1)$ and $f^\ve\in\cH$ one has
 \begin{equation*}
  \frac{c_\eta}{2}\,\|f^\ve\|_{L^2(\dD x\dD \gamma)}^2\,\leq\, H(f^\ve)\,\leq \,\frac{C_\eta}{2}\,\|f^\ve\|_{L^2(\dD x\dD \gamma)}^2\,.
 \end{equation*}
 \label{lem:equivnorm}
\end{lem}
\begin{proof}
 The result follows directly from the expression of the entropy and the estimate
 \[
  |\lla j^\ve, \partial_x\phi^\ve\rra_{L^2_x}|\,\leq\,\|j^\ve\|_{L^2_x}\,\|\partial_x\phi^\ve\|_{L^2_x}\,\leq\,\frac{1}{\ve}\,C_P\,m_2^{1/2}\,\|f^\ve\|^2_{L^2(\dD x\dD \gamma)}\,.
 \]
\end{proof}

With the previous lemmas, we can prove the main result of this section.

\begin{proof}[Proof of Proposition~\ref{prop:hypocoercivity}]
Using the moment equation $\ve^2\partial_tj^\ve+ m_2\,\partial_x\rho^\ve + \partial_xS^\ve  + j^\ve = 0$, the time derivative of the entropy splits into five terms
\[
\frac{\dD}{\dD t}H(f^\ve(t))\ =\ T_1^\ve(t) + \eta\,m_2\,T_2^\ve(t) + \eta\,T_3^\ve(t) + \eta\,T_4^\ve(t) + \eta\,T_5^\ve(t)\,.
\]
Let us precise and estimate each term:
\begin{eqnarray*}
 T_1^\ve\,:=\,\frac{1}{2}\frac{\dD}{\dD t}\|f^\ve\|^2_{L^2(\dD x\dD \gamma)}&\leq&-\frac{1}{\ve^2}\|f^\ve - \rho^\ve M\|_{L^2(\dD x\dD \gamma)}^2\,,\\
T_2^\ve\,:=\,-\lla\partial_x\rho^\ve,\partial_x\phi^\ve\rra_{L^2_x} = \lla\rho^\ve,\partial_{xx}^2\phi^\ve\rra_{L^2_x}&=&-\|\rho^\ve\|^2_{L^2_x}\,,\\
T_3^\ve\,:=\,-\lla\partial_xS^\ve,\partial_x\phi^\ve\rra\,=\,-\lla S^\ve,\rho^\ve\rra&\leq&(m_4 - m_2^2)^{1/2}\,\|f^\ve - \rho^\ve M\|_{L^2(\dD x\dD \gamma)}\,\|\rho^\ve\|_{L^2_x}\,,\\
T_4^\ve\,:=\,-\lla j^\ve,\partial_x\phi^\ve\rra&\leq&\,\frac{C_P\,m_2^{1/2}}{\ve}\,\|f^\ve - \rho^\ve M\|_{L^2(\dD x\dD \gamma)}\,\|\rho^\ve\|_{L^2_x}\,,\\
T_5^\ve\,:=\,\ve^2\lla j^\ve,\partial_t\partial_x\phi^\ve\rra&\leq&m_2\|f^\ve - \rho^\ve M\|_{L^2(\dD x\dD \gamma)}^2\,,
\end{eqnarray*}
where we used the estimates of Lemma~\ref{lem:moments} and Lemma~\ref{lem:estimphi}. Combining everything, one has 
\begin{multline*}
 \frac{\dD}{\dD t}H(f^\ve(t))\ +(\frac{1}{\ve^2}-\eta\,m_2)\,\|f^\ve - \rho^\ve M\|_{L^2(\dD x\dD \gamma)}^2\,+\,\eta\,\,m_2\|\rho^\ve\|^2_{L^2_x}\\
 \leq\ \eta\,\left((m_4 - m_2^2)^{1/2} + \frac{C_P\,m_2^{1/2}}{\ve}\right)\|f^\ve - \rho^\ve M\|_{L^2(\dD x\dD \gamma)}\,\|\rho^\ve\|_{L^2_x}\,.
\end{multline*}
It follows that for any $\eta\in(0,\eta_2)$ with $\eta_2 = m_2/(((m_4-m_2^2)^{1/2}+C_P\,m_2^{1/2})^2 + m_2^2)$ and $\ve\in(0,1)$, one has %there is a constant $\kappa_\eta$ such that 
\[
 \frac{\dD}{\dD t}H(f^\ve(t)) + K_\eta\,(\|f^\ve - \rho^\ve M\|_{L^2(\dD x\dD \gamma)}^2 + \|\rho^\ve M\|^2_{L^2(\dD x\dD \gamma)})\,\leq\,0\,,
\]
with $K_\eta = \min(1-\eta m_2,\eta m_2)/2$. Then use the triangle inequality,  $(a+b)^2\leq 2(a^2+b^2)$ and Lemma~\ref{lem:equivnorm} to get that
\[
\frac{\dD}{\dD t}H(f^\ve(t)) + \frac{K_\eta}{C_\eta}\,H(f^\ve(t))\,\leq\,0\,,
\]
for all $\eta\in(0,\min(\eta_1,\eta_2))$. Choose some admissible $\eta$ and set $\kappa = K_\eta/(2C_\eta)$ to conclude.
\end{proof}

\section{The fully discrete setting}
\label{sec:FullDiscrete}

In this section, our aim is now to adapt the previous results to the discrete setting. To proceed, we first introduce an implicit in time and finite volume in space discretization of \eqref{eq:Collision} both for Fokker-Planck and BGK cases. We then prove the discrete counterparts of the results presented in the previous section, namely we study the diffusive limit and the exponential return to equilibrium.

\subsection{Notations}~

\mysubsection{Mesh.} We first restrict the velocity domain to a bounded symmetric segment $[-v_{\star},v_{\star}]$, since it is not possible in practice to implement a numerical scheme on an unbounded domain. We consider a primal mesh of this interval composed of $2L$ control volumes arranged symmetrically around $v = 0$. We thus get $2L+1$ distinct interface points denoted by $v_{j+\frac{1}{2}}$ (for consistency with usual finite volume notations) with $j = -L,\dots, L$. In this way
\[
v_{-L+1/2} = -v_{\star}\,,\quad v_{1/2} = 0\,,\quad v_{j+1/2} = -v_{-j+1/2}\quad \forall j= 0,\dots, L\,.
\]
The cells of the primal mesh are given by 
\[
\V_j\ :=\ (v_{j-\frac{1}{2}},v_{j+\frac{1}{2}})\,,\quad j\in\J \,:=\, \{-L+1,\dots,L\}.
\]
Each cell $\V_j$ has length $\Delta v_{j}=v_{j+\frac{1}{2}}-v_{j-\frac{1}{2}}$ and midpoint $v_j$. At the level of cells, the symmetry of the velocities reads $v_{j}\ =\ - v_{-j+1}$ for all $j = 1,\dots,L$. We also define the dual mesh with cells  
\[
\V_{j+1/2}^*\ :=\ (v_{j},v_{j+1})\,,\quad j\in\J^* \,:=\, \{-L,\dots,L\}\,,
\]
with $v_{-L} := v_{-L+1/2} = -v_{\star}$ and $v_{L+1} := v_{L+1/2} = v_{\star}$. The length of the dual cell $\V_{j+1/2}^*$ is $\Delta v_{j+1/2} = v_{j+1} - v_{j}$. Notations introduced are illustrated in Fig. \ref{fig:domaine_v}.

\begin{figure}
\centering
\begin{tikzpicture}[scale=1.]
\draw[line width =0.7pt, color=blue] (6,-6) -- (10,-6);
\draw[line width =0.7pt, color=blue] (3,-6) -- (4,-6);
\draw[line width =0.7pt, color=blue] (12,-6) -- (13,-6);
\draw[line width =0.7pt,densely dashed,color=blue] (4,-6) -- (6,-6);
\draw[line width =0.7pt,densely dashed,color=blue] (10,-6) -- (12,-6);
\path (4,-6)  node{\scriptsize\color{red} $\blacklozenge$};
\path (4,-6.25)  node{\small\color{red}$v_{-L+1}$};
\path (8,-6)  node{\scriptsize $\bullet$};
\path (8,-5.8)  node{\small $0$};
\path (13,-5.8)  node{\small $v_*$};
\path (3,-5.8)  node{\small $-v_*$};
\path (2.8,-6.3)  node{\small $v_{-L+1/2}$};
\path (3,-6.7)  node{\small\color{red} $=v_{-L}$};
\path (8,-6.3)  node{\small $v_{1/2}$};
\path (7.375,-6)  node{\scriptsize\color{red} $\blacklozenge$};
\path (7.375,-6.25)  node{\small\color{red} $v_{0}$};
\path (8.625,-6)  node{\scriptsize\color{red} $\blacklozenge$};
\path (8.625,-6.25)  node{\small\color{red} $v_{1}$};
\path (12,-6)  node{\scriptsize\color{red} $\blacklozenge$};
\path (13,-6)  node{\scriptsize $\bullet$};
\path (3,-6)  node{\scriptsize $\bullet$};
\path (13,-6.3)  node{\small $v_{L+1/2}$};
\path (13,-6.7) node{\small \color{red} $=v_{L+1}$};
\path (6.75,-6)  node{\scriptsize $\bullet$};
\path (6.75,-6.3)  node{\small $v_{-1/2}$};
\path (9.25,-6)  node{\scriptsize $\bullet$};
\path (9.25,-6.3)  node{\small $v_{3/2}$};
\path (12,-6.25)  node{\small\color{red}$v_{L}$};
\path (6,-6) node{\scriptsize\color{red} $\blacklozenge$};
\path (10,-6) node{\scriptsize\color{red} $\blacklozenge$};
\path (6,-6.25)  node{\small\color{red} $v_{-1}$};
\path (10,-6.25)  node{\small\color{red} $v_{2}$};
\draw[<->,line width =0.9pt, densely dashed,color=black](7.375,-6.5)--(8.625,-6.5) ;
\draw[<->,line width =0.9pt, densely dashed,color=black](8,-5.5)--(9.25,-5.5) ;
\path (8.625,-5.3)  node{\small\color{black} $\Delta v_{1}$};
\path (8,-6.8)  node{\small\color{black} $\Delta v_{1/2}$};
\end{tikzpicture}
\caption{Discretization of the velocity domain.}
\label{fig:domaine_v}
\end{figure}

In space, we consider a discretization of the torus $\T$ into $N$ subintervals 
\[
\X_i\ :=\ (x_{i-\frac{1}{2}},x_{i+\frac{1}{2}})\,,\quad i\in\I\,:=\,\Z/N\Z
\]
of length $\Delta x_{i}$ and centers $x_{i}$. In what follows, we assume that $N$ is odd. Indeed, this assumption is natural to obtain, among others, a discrete Poincaré inequality on the torus with our choice of discrete gradients, as explained later.

The control volumes in phase space are defined by 
\[
 K_{ij}\ :=\ \X_i\times\V_j\,,\qquad\forall (i,j)\in\I\times\J.
\]
The size of the phase-space discretization is defined by $\delta\ =\ \max(\Delta x, \Delta v)$ where $\Delta x$ and $\Delta v$ are the maximum of $(\Delta x_i)_{i\in\I}$ and $(\Delta v_j)_{j\in\J}$ respectively. Finally, we set $\Delta t>0$ the time step, and we define $t^{n}=n\Delta t$ for all $n\geq 0$.

\mysubsection{Discrete Maxwellians.} In the case of the BGK operator \eqref{def:Qbgk}, we assume that we are given cell values $\M=(\M_j)_{j\in\J}\in\RR^{\J}$ satisfying the following assumptions
\begin{equation}\label{hyp:MaxwBGK}
\left\{
\begin{array}{l}
\M_j > 0\,,\quad\M_{j}\ =\ \M_{-j+1}\,,\quad\forall j = 1,\dots,L\,;\\[.75em]
\ds\sum_{j\in\J}\M_j\,\Delta v_j\ =\ 1\,;\\[.75em]
\ds0\,<\,\underline{m}_2\,\leq\, m_2^{\Delta v}\,\leq\,\overline{m}_2\,,\quad m_4^{\Delta v}\,\leq\, \overline{m}_4\,,
\end{array}
\right.
\end{equation}
where for $k\in\NN$
\[
 m_k^{\Delta v}\ :=\ \sum_{j\in\J}|v_j|^k\,\M_j\,\Delta v_j
\]
and $\underline{m}_2$, $\overline{m}_2$, $\overline{m}_4$ are some universal constants.

In the case of the Fokker-Planck operator \eqref{def:Qfp}, we assume that we are given interface values $(\M_{j+1/2}^*)_{j\in\J^*}\in\RR^{\J^*}$ such that
\begin{equation}\label{hyp:MaxwFP}
\left\{
\begin{array}{l}
\M_{j+1/2}^*\ =\ \M_{-j+1/2}^*\,,\quad\forall j \in\J^*\,;\\[.75em]
\M_{L+1/2}^*\ =\ \M_{-L+1/2}^*\ =\ 0\,;\\[.75em]
\ds\M_j\ :=\ \frac{\M_{j-1/2}^* - \M_{j+1/2}^*}{v_j\,\Delta v_j}\,>\,0\,, \quad\forall j\in\J\,;\\[.75em]
\ds\sum_{j\in\J}\M_j\,\Delta v_j\ =\ 1\,;\quad\\[.75em]
\ds0\,<\,\underline{m}_2\,\leq\, m_2^{\Delta v}\,\leq\,\overline{m}_2\,,\quad m_4^{\Delta v}\,\leq\, \overline{m}_4\,.
\end{array}
\right.
\end{equation}

These assumptions on the discrete Maxwellians $\M$ are weak and it is easy to build an example.
\begin{example} \label{example_M} Define $\M_j := c_{\Delta v}\,\M(v_j)$ in the BGK case and $\M_{j+1/2} := \tilde{c}_{\Delta v}\,\M(v_{j+1/2})$ in the Fokker-Planck case and compute $c_{\Delta v}$ and $\tilde{c}_{\Delta v}$ to normalize the mass of $(\M_j)_j$ to $1$. Then, by consistency of the piecewise constant approximation (say, in $L^1$) the last condition in \eqref{hyp:MaxwBGK} (\textit{resp.} \eqref{hyp:MaxwFP}) holds in the BGK case (\textit{resp.} the Fokker-Planck case).
\end{example}

\begin{rem}
  ~ 
  \begin{itemize}
  \item The second assumption in \eqref{hyp:MaxwFP} is needed for technical reasons, in order to perform discrete integration by parts. 
   However, one can note that for sufficiently large domains this hypothesis is relevant and invisible in practice because of the fast decay of the Gaussian.
    \item In both cases observe that symmetries of both the velocities and the Maxwellians imply
      \begin{equation}\label{zero_mean_M}
        \sum_{j\in\J}\Delta v_{j}v_{j}\M_{j}=0.
      \end{equation}

    \item The third assumption in  \eqref{hyp:MaxwFP} is the discrete counterpart of \eqref{eq:MaxwellianIdentity}, namely $\pa_{v}\M=-v\M$.
    \item Let us underline that $m_2^{\Delta v}$ and $m_4^{\Delta v}$ are not consistent with $m_2$ and $m_4$ when $\Delta v$ tends to zero, since we consider the problem in a bounded domain $[-v_*,v_*]$ with zero boundary values for the discrete Maxwellian. We then have an error in $O(\M(v_*))$ for these quantities.
  \end{itemize}
\end{rem}
As in the continuous case, we also introduce the inverse of the Maxwellian 
\[
\gamma_j\ =\ \frac{1}{\M_j}\,,\quad \forall j\in\J\,.
\]

\mysubsection{Discrete gradients.} Given some discrete microscopic quantity $g = (g_{ij})_{i\in\I,j\in\J}$ defined on control volumes and assuming that some boundary data (in velocity) $(g_{i, L+1})_{i\in\I}$ and $(g_{i, -L})_{i\in\I}$ are given, one may define the discrete gradients $D_v g \in\RR^{\I\times\J^*}$ on the dual mesh by
\[
 (D_v g)_{i,j+1/2}\ =\ \frac{g_{i,j+1} - g_{ij}}{\Delta v_{j+1/2}}\,,\quad \forall j\in\J^*\,,\quad\forall i \in\I.
\]
In the numerical analysis of our scheme it will be convenient to use discrete gradients in space too. Given a macroscopic discrete quantity $\rho = (\rho_{i})_{i\in\I}$, we define discrete gradients $D_x\rho \in\RR^{\I}$ on the (primal) mesh as
\[
 (D_x \rho)_i\ =\ \frac{\rho_{i+1/2} - \rho_{i-1/2}}{\Delta x_{i}}\,,\quad \text{where}\quad \rho_{i+1/2}\ =\ \frac{\rho_{i+1}+\rho_{i}}{2}\,,\quad \forall i\in\I\,.
\]

\mysubsection{Discrete functional spaces.} From there we define for microscopic quantities the discrete weighted $L^2$ norm
\[
 \|g\|_{2,\gamma}^2\ :=\
 \sum_{(i,j)\in\I\times\J}|g_{ij}|^2\,\gamma_{j}\,\Delta x_i\,\Delta v_j\,.
\]
For macroscopic quantities, we shall need the discrete $L^2$ norm
\[
 \|\rho\|_2^2\ :=\  
 \sum_{i\in\I}|\rho_{i}|^2\,\Delta x_i\,
\]
and discrete Sobolev seminorm
\[
 \|D_x\rho\|_2^2\ :=\ 
 \sum_{i\in\I}|(D_x\rho)_{i}|^2\,\Delta x_i\,.
\]
Finally, for a discrete Maxwellian satisfying \eqref{hyp:MaxwFP} one can introduce the Sobolev seminorm in velocity
\[
\|D_vg\|_{2,\M^*}^2\ :=\ 
\sum_{(i,j)\in\I\times\J^*}\left|(D_vg)_{i,j+1/2}\right|^2\,\M_{j+1/2}^*\,\Delta x_i\,\Delta v_{j+1/2}\,.
\]
Observe that since the discrete Maxwellian vanishes at endpoints of the domain, the definition does not involve boundary data  $(g_{i, L+1})_{i\in\I}$ and $(g_{i, -L})_{i\in\I}$.  One has the following discrete counterparts of the Gaussian and flat Poincaré inequalities. The proof of these results is given in appendix. 
\begin{lem}[Discrete Gaussian Poincaré inequality on bounded velocity domain]\label{lem_poincare_v}
Given a discrete Maxwellian on interfaces $(\M_{j+1/2}^*)_{j\in\J^*}$ satisfying \eqref{hyp:MaxwFP}, one has for all $f = (f_{j})_{j\in\J}\in\RR^\J$ that
\begin{equation}\label{ineq_Poincare_v}
\|f - \rho\,\M \|_{2,\gamma}^2\ \leq\ \left\|D_v\left(\frac{f}{\M}\right)\right\|_{2,\M^*}^2,
\end{equation}
where $\M = (\M_{j})_{j\in\J}$ is defined in \eqref{hyp:MaxwFP} and $\rho = \sum_{j\in\J} f_j \Delta v_j$.
\end{lem}

\begin{lem}[Discrete Poincaré inequality on the torus]
Assume that the number of points $N$ in the space discretization of the torus is odd. Then, there is a constant $C_P$ that does not depend on $\Delta x$ such that one has for all $\phi = (\phi_{i})_{i\in\I}\in\RR^\I$ satisfying $\sum_{i\in \I}\Delta x_i\phi_i\ =\ 0$ that 
\begin{equation}\label{ineq_Poincare_x}
\|\phi\|_{2}\, \leq\,C_P\, \|D_x\phi\|_{2}.
\end{equation}
\label{lem:poincare_discrete}
\end{lem}

\subsection{Presentation of the schemes}

In this section we present numerical schemes for Equation~\eqref{eq:Collision} in  both the Fokker-Planck and the BGK cases. The schemes are of finite volume type \cite{eymard_2000_finite}, meaning that they are based on an approximation of the fluxes appearing in the integrated version of \eqref{eq:Collision} on each cell $K_{ij}$. In the time variable we choose a backward Euler discretization. First of all, we discretize the initial datum $f_{0}$ by
\begin{equation*}
f_{ij}^{0}=\frac{1}{\Delta x_{i}\,\Delta v_{j}}\iint_{K_{ij}}f(0,x,v)\,\dD x\,\dD v,\quad \forall (i,j)\in\I\times\J.
\end{equation*}
Now let us define the schemes corresponding to each collision operator. To lighten the notations, we omit the superscript $\varepsilon$ for the discretization of $f^\ve$ and its moments in all this section, even if these quantities obviously depend on $\varepsilon$ as in the continuous setting.

\mysubsection{Fully discrete Fokker-Planck equation.} In the case of the Fokker-Planck operator \eqref{def:Qfp_grad}, the scheme is given as follows. For all $i\in\I$, $j\in\J$, $n\geq 0$,
\begin{equation}\label{eq:scheme_f_FP}	
\varepsilon\Delta x_{i}\Delta v_{j}(f_{ij}^{n+1}-f_{ij}^{n})+\Delta t\left({\F}_{i+\frac{1}{2},j}^{n+1}-{\F}_{i-\frac{1}{2},j}^{n+1}\right)=\frac{\Delta t}{\varepsilon}\left({\G}_{i,j+\frac{1}{2}}^{n+1}-{\G}_{i,j-\frac{1}{2}}^{n+1}\right),
\end{equation}
where the numerical fluxes are respectively defined for all $n\geq 0$ by
\begin{align}
{\F}_{i+\frac{1}{2},j}^{n+1}\ =\ v_{j}\frac{f_{i+1,j}^{n+1}+f_{ij}^{n+1}}{2}\Delta v_{j},&&&\forall j\in\J,\quad\forall i\in\I\,,\label{eq:def_fluxF_f}\\[1em]
{\G}_{i,j+\frac{1}{2}}^{n+1}\ =\ \M_{j+\frac{1}{2}}^*\left(\frac{f_{i,j+1}^{n+1}}{\M_{j+1}}-\frac{f_{ij}^{n+1}}{\M_{j}}\right)\frac{1}{\Delta v_{j+\frac{1}{2}}}\Delta x_{i},&&& \forall j\in\J^*\setminus\{-L,L\},\quad\forall i\in\I\,,\label{eq:def_fluxG_f}\\[1em]
{\G}_{i,-L+\frac{1}{2}}^{n+1}\ =\ {\G}_{i,L+\frac{1}{2}}^{n+1}\ =\ 0,&&&\forall i\in\I\,.\label{eq:fluxnul_f}
\end{align}

\begin{rem}
  The first equation \eqref{eq:def_fluxF_f} is a centered discretization of the free transport. The choice of the velocity fluxes $\G_{i,j+\frac{1}{2}}$ in \eqref{eq:def_fluxG_f}, which is based on the gradient form \eqref{def:Qfp_grad} of the Fokker-Planck operator
  \[
  \Q_{FP}(f)\ =\ \partial_v\left( \M\partial_v\left(f/\M\right)\right) 
  \]
  is fairly close to the Chang-Cooper \cite{ChangCooper:1970, BuetDellacherie:2010} or Scharfetter-Gummel \cite{IlIn1969,Scharfetter1969} approximation. Indeed, if we consider a uniform velocity discretization with $\Delta v_j=\Delta v$ constant, and the discrete Maxwellian defined as in Example \ref{example_M}, the velocity flux is then given by
  \[
  \G_{i,j+\frac{1}{2}}=\frac{\Delta x_i}{\Delta v}\left(B\left(-v_{j+1}\Delta v\right)f_{i,j+1}-B\left(v_j\Delta v\right)f_{i,j}  \right),
  \]
  where $B$ is the Bernoulli function defined by $B(x)=x/(e^x-1)$ if $x\neq 0$, $B(0)=1$. The Chang-Cooper and Scharfetter-Gummel schemes are obtained by replacing $v_{j}$ and $v_{j+1}$ by $v_{j+1/2}$.
  
  Finally, the last equation \eqref{eq:fluxnul_f} corresponds to zero flux boundary conditions in velocity.
  
\end{rem}

\mysubsection{Fully discrete BGK equation.} In the case of the BGK operator \eqref{def:Qbgk}, the scheme is given as follows. For all $i\in\I$, $j\in\J$, $n\geq 0$,
\begin{equation}\label{eq:scheme_f_BGK}	
\varepsilon\Delta x_{i}\Delta v_{j}(f_{ij}^{n+1}-f_{ij}^{n})+\Delta t\left({\F}_{i+\frac{1}{2},j}^{n+1}-{\F}_{i-\frac{1}{2},j}^{n+1}\right)=\frac{\Delta t}{\varepsilon}\Delta x_{i}\Delta v_{j}\left(\rho_{i}^{n+1}\M_{j}-f_{ij}^{n+1}\right),
\end{equation}
where the numerical flux ${\F}_{i+\frac{1}{2},j}^{n+1}$ is still defined by \eqref{eq:def_fluxF_f} and for all $i\in\I$ and $n\geq 0$, the discrete macroscopic density is given by
\begin{equation}\label{eq:def_rho}
\rho_{i}^{n}\ :=\ \sum_{j\in\J}\Delta v_{j}f_{ij}^{n}\,.
\end{equation}
For future use and independently of the choice of collision operator, we also define the other used velocity moments of the discrete distribution, namely
\begin{equation}\label{eq:def_jS}
 J_{i}^{n}\ :=\ \frac{1}{\ve}\,\sum_{j\in\J}\Delta v_{j}\,v_j\,f_{ij}^{n}\,,\qquad S_{i}^{n}\ :=\ \sum_{j\in\J}\Delta v_{j}\,(v_j^2-m_2^{\Delta v})\,f_{ij}^{n}\,,
\end{equation}
for all $i\in\I$ and $n\geq 0$.
Both schemes \eqref{eq:scheme_f_FP} and \eqref{eq:scheme_f_BGK} clearly satisfy the discrete mass conservation:
\begin{equation*}
    \sum_{(i,j)\in\I\times\J}\Delta x_i\Delta v_j f_{i,j}^n= \sum_{(i,j)\in\I\times\J}\Delta x_i\Delta v_j f_{i,j}^0 = \mu_f, \quad\forall n\geq 0.
\end{equation*}
This property is obtained by using the zero flux boundary conditions in velocity \eqref{eq:fluxnul_f} in the Fokker-Planck case, and assumption on the discrete Maxwellian \eqref{hyp:MaxwBGK} in the BGK case. 

Existence and uniqueness of a solution to the fully implicit schemes \eqref{eq:scheme_f_FP} and \eqref{eq:scheme_f_BGK} will be established in the next section as a by-product of the discrete ``entropy'' estimate given in Lemma~\ref{lem:discrete_entropy}.

\subsection{Uniform estimates}

In this section we derive the discrete counterparts of the results of Section~\ref{sec:unif}.

\begin{lem}[Discrete ``entropy'' estimate]\label{lem:discrete_entropy}
 Let us consider a discrete Maxwellian satisfying \eqref{hyp:MaxwFP} (\text{resp.} \eqref{hyp:MaxwBGK}) and let $(f_{ij}^n)_{i\in\I,j\in\J,n\in\NN}$ solve the scheme \eqref{eq:scheme_f_FP} (\text{resp.} \eqref{eq:scheme_f_BGK}). Then for every $n\geq 0$\,,

 \begin{equation}
 \frac{\|f^{n+1}\|^2_{2,\gamma} - \|f^n\|^2_{2,\gamma}}{2\Delta t}\,+\,\frac{1}{\ve^2}\|f^{n+1}-\rho^{n+1} \M\|^2_{2,\gamma}\ \leq\ 0\,.
  \label{eq:estimEI_dis}
 \end{equation}
 In particular one has
 \begin{equation}
  \max\left(\sup_{n\geq0}\|f^n\|^2_{2,\gamma},\ \frac{2}{\ve^2}\sum_{n=1}^\infty\Delta t\,\|f^{n}-\rho^{n} \M\|^2_{2,\gamma}\right)\,\leq\,\|f^0\|^2_{2,\gamma}\,.
  \label{eq:L2gamma_discrete}
 \end{equation}

\label{lem:basic_estimates_discrete}
\end{lem}

Since schemes \eqref{eq:scheme_f_FP} and \eqref{eq:scheme_f_BGK} are finite dimensional linear systems, we deduce from \eqref{eq:L2gamma_discrete} the uniqueness of a solution, and then the existence.

\begin{cor}
The scheme \eqref{eq:scheme_f_FP} (resp. \eqref{eq:scheme_f_BGK}) admits a unique solution $(f_{ij}^n)_{i\in\I,j\in\J,n\in\NN}$.
\end{cor}

\begin{proof}[Proof of Lemma \ref{lem:discrete_entropy}]
We first consider the Fokker-Planck case. Multiplying \eqref{eq:scheme_f_FP} by $f_{ij}^{n+1}\gamma_{j}$ and summing over $i$ and $j$, we get
\begin{gather}\label{eq:estimEI1_FP}
\varepsilon\sum_{(i,j)\in\I\times\J}\Delta x_{i}\Delta v_{j}(f^{n+1}_{ij}-f_{ij}^{n})f_{ij}^{n+1}\gamma_{j}+\Delta t\sum_{(i,j)\in\I\times\J}(\F_{i+\frac{1}{2},j}^{n+1}-\F_{i-\frac{1}{2},j}^{n+1})f_{ij}^{n+1}\gamma_{j}=\\
\frac{\Delta t}{\varepsilon}\sum_{(i,j)\in\I\times\J}(\G_{i,j+\frac{1}{2}}^{n+1}-\G_{i,j-\frac{1}{2}}^{n+1})f_{ij}^{n+1}\gamma_{j}.\nonumber
\end{gather}
Using that $a(a-b)\geq (a^2-b^2)/2$ for all $a$, $b\in\R$, the first term of \eqref{eq:estimEI1_FP} can be bounded from below by
\[
 \frac{\varepsilon}{2}\sum_{(i,j)\in\I\times\J}\Delta x_{i}\Delta v_{j}\left((f_{ij}^{n+1})^2-(f_{ij}^n)^2\right)\gamma_{j}=\frac{\varepsilon}{2}(\|f^{n+1}\|^2_{2,\gamma} - \|f^n\|^2_{2,\gamma}).
\]
Then, using the definition \eqref{eq:def_fluxF_f} of the numerical fluxes, the second term of \eqref{eq:estimEI1_FP} gives
\begin{align*}
\sum_{(i,j)\in\I\times\J}(\F_{i+\frac{1}{2},j}^{n+1}-\F_{i-\frac{1}{2},j}^{n+1})f_{ij}^{n+1}\gamma_j&=\sum_{(i,j)\in\I\times\J}\Delta v_{j}v_{j}\gamma_{j}\left(\frac{f_{i+1,j}^{n+1}+f_{ij}^{n+1}}{2}-\frac{f_{ij}^{n+1}+f_{i-1,j}^{n+1}}{2}\right)f_{ij}^{n+1}\\
&=\sum_{j\in\J}\Delta v_{j}v_{j}\gamma_{j}\sum_{i\in\I}\frac{f_{i+1,j}^{n+1}f_{ij}^{n+1}-f_{ij}^{n+1}f_{i-1,j}^{n+1}}{2}=0.
\end{align*}
Finally, we perform a discrete integration by parts in velocity on the third term of \eqref{eq:estimEI1_FP} using the zero flux boundary conditions, to obtain
\[
\begin{array}{rl}
\ds\sum_{(i,j)\in\I\times\J}(\G_{i,j+\frac{1}{2}}^{n+1}-\G_{i,j-\frac{1}{2}}^{n+1})f_{ij}^{n+1}\gamma_j&=\,\ds-\sum_{(i,j)\in\I\times\J^*}\G_{i,j+\frac{1}{2}}^{n+1}(f_{i,j+1}^{n+1}\gamma_{j+1}-f_{ij}^{n+1}\gamma_j)\\[1em]
&=\,\ds-\sum_{(i,j)\in\I\times\J^*}\frac{\Delta x_{i}}{\Delta v_{j+\frac{1}{2}}}\M_{j+\frac{1}{2}}^*(f_{i,j+1}^{n+1}\gamma_{j+1}-f_{ij}^{n+1}\gamma_j)^2\\[1.5em]
&=\,\ds- \|D_v(f/\M)\|^2_{2,\M^*}\\[1em]
&\ds\leq\,- \|f^{n+1}-\rho^{n+1} \M\|^2_{2,\gamma}\,,
\end{array}
\]
thanks to the discrete Gaussian Poincaré inequality. This concludes the proof in the Fokker-Planck case.

For the BGK operator the left hand side is exactly the same as in \eqref{eq:estimEI1_FP}. Using the definition \eqref{eq:def_rho} of $\rho_{i}^{n+1}$ and the unit mass assumption on the discrete Mawxellian in \eqref{hyp:MaxwBGK}, we obtain for the right hand side
\begin{align*}
\sum_{(i,j)\in\I\times\J}\Delta x_{i}\Delta v_{j}(\rho_{i}^{n+1}\M_j-f_{ij}^{n+1})f_{ij}^{n+1}\gamma_{j}&=-\sum_{(i,j)\in\I\times\J}\Delta x_{i}\Delta v_{j}(f_{ij}^{n+1}-\rho_{i}^{n+1}\M_j)^2\gamma_{j}\\
&=- \|f^{n+1}-\rho^{n+1} \M\|^2_{2,\gamma},
\end{align*}
which concludes the proof.
\end{proof}

\begin{lem}[Discrete moments estimates]
 Under the hypotheses of Lemma~\ref{lem:basic_estimates_discrete}, the moments $(\rho_i^n)_{i\in\I}$, $(J_i^n)_{i\in\I}$ and $(S_i^n)_{i\in\I}$ satisfy  the following estimates, for all $n\in\NN$:
 \begin{equation}
  \|\rho^n\|_{2}\ \leq\ \|f^n\|_{2,\gamma}\,,
  \label{eq:estmomentsLinfrho_discrete}
 \end{equation}
 \begin{equation}
 \ve\,\|J^n\|_{2}\ \leq\ (m_2^{\Delta v})^{1/2}\,\|f^n\|_{2,\gamma}\,,
  \label{eq:estmomentsLinfj_discrete}
 \end{equation}
 \begin{equation}
 \ve\,\|J^n\|_{2}\ \leq\ (m_2^{\Delta v})^{1/2}\,\|f^n-\rho^n\,\M\|_{2,\gamma}\,,
  \label{eq:estmomentsL2j_discrete}
 \end{equation}
 \begin{equation}
 \|S^n\|_{2}\ \leq\ (m_4^{\Delta v} - (m_2^{\Delta v})^2)^{1/2}\,\|f^n-\rho^n\,\M\|_{2,\gamma}\,.
  \label{eq:estmomentsL2S_discrete}
 \end{equation}
 In particular, there is a constant $C>0$ depending only on the uniform bounds on $m_2^{\Delta v}$ and $m_4^{\Delta v}$ from hypotheses \eqref{hyp:MaxwBGK} and \eqref{hyp:MaxwFP} such that
 \begin{multline}
\max\left(\sup_{n\in\NN}\|\rho^n\|_{2},\,\ve\,\sup_{n\in\NN}\|J^n\|_2,\,\left(\sum_{n=0}^\infty\Delta t\,\|J^n\|_2^2\right)^{1/2},%\right.\\\left.
\, \frac1\ve\,\left(\sum_{n=0}^\infty\Delta t\,\|S^n\|_2^2\right)^{1/2}\right)\\\ \leq\ C\,\|f^0\|_{2,\gamma}\,.
\label{eq:estimmoments_discrete}
 \end{multline}
\end{lem}
\begin{proof}
The proof is almost identical to that of Lemma~\ref{lem:moments}. Indeed, just observe that the moments rewrite as follows. First one has
\[
 \rho_i^n\ =\ \sum_{j\in\J}\Delta v_j\,(f^n_{ij}\gamma_j^{1/2})\M_j^{1/2}\,.
\]
Then, thanks to the symmetry of the discrete Maxwellian and anti-symmetry of discrete velocities yielding \eqref{zero_mean_M}, one has 
\[
 J_i^n\ =\ \frac{1}{\ve}\sum_{j\in\J}\Delta v_j\,(f^n_{ij}\gamma_j^{1/2})\,(v_j\,\M_j^{1/2})\ =\ \frac{1}{\ve}\sum_{j\in\J}\Delta v_j\,(f^n_{ij}-\rho_i^n\,\M_j)\gamma_j^{1/2}\,(v_j\,\M_j^{1/2})\,.
\]
Finally,  observe that by definition of $m_2^{\Delta v}$ one has
\[
 S_i^n\ =\ \sum_{j\in\J}\Delta v_j\,(f^n_{ij}-\rho_i^n\,\M_j)\gamma_j^{1/2}\,(v_j^2-m_2^{\Delta v})\,\M_j^{1/2}\,.
\]
Then one can take the norm of each expression and use Cauchy-Schwarz inequalities to get the desired estimates. Estimate \eqref{eq:estimmoments_discrete} is then a direct consequence of \eqref{eq:L2gamma_discrete}. 
\end{proof}

\subsection{Asymptotic-preserving property}

For macroscopic quantities we recall that approximate values at the interfaces $x_{i+\frac{1}{2}}$, $i\in\I$ are defined by the average
\[Q_{i+\frac{1}{2}}=\frac{Q_{i}+Q_{i+1}}{2},\quad X=\rho,\,J,\,S.\]

\begin{lem}[Moments equations]
Let us consider a discrete Maxwellian satisfying \eqref{hyp:MaxwFP} (\text{resp.} \eqref{hyp:MaxwBGK}) and a solution to the scheme \eqref{eq:scheme_f_FP} (\text{resp.} \eqref{eq:scheme_f_BGK}). Then the discrete moments satisfy the following equations. For all $i\in\I$, $n\geq 0$,
\begin{align}
\label{scheme_rho}\Delta x_{i}\,(\rho_{i}^{n+1}-\rho_{i}^{n})\,+\,\Delta t\,(J_{i+\frac{1}{2}}^{n+1}-J_{i-\frac{1}{2}}^{n+1})&& =&&& 0,\\[.5em]
\label{scheme_J}\varepsilon^2\,\Delta x_{i}\,(J_{i}^{n+1}-J_{i}^n)\,+\,\Delta t\,(S_{i+\frac{1}{2}}^{n+1}-S_{i-\frac{1}{2}}^{n+1})\,+\,\Delta t\,m_{2}^{\Delta v}\,(\rho_{i+\frac{1}{2}}^{n+1}-\rho_{i-\frac{1}{2}}^{n+1})&& =&&&-\Delta t\,\Delta x_{i}\,J_{i}^{n+1}.
\end{align}

\end{lem}

\begin{proof}
Equation \eqref{scheme_rho} can be easily obtained by taking \eqref{eq:scheme_f_FP} (resp. \eqref{eq:scheme_f_BGK}) and summing over $j\in\J$. Thanks to the zero flux boundary conditions in the Fokker-Planck case, and the  normalization assumption of the Maxwellian in the BGK case, one easily sees that the contribution of the collision operators vanishes.

To obtain the second equation, we multiply \eqref{eq:scheme_f_FP} (resp. \eqref{eq:scheme_f_BGK}) by $v_{j}$ and sum over $j$. In the Fokker-Planck case, it yields
\begin{gather*}
\varepsilon^2\Delta x_{i}(J_{i}^{n+1}-J_{i}^{n})+\Delta t(S_{i+\frac{1}{2}}^{n+1}-S_{i-\frac{1}{2}}^{n+1})+\Delta tm_{2}^{\Delta v}(\rho_{i+\frac{1}{2}}^{n+1}-\rho_{i-\frac{1}{2}}^{n+1})=\\
\frac{\Delta t}{\varepsilon}\sum_{j\in\J}v_{j}\left(\G_{i,j+\frac{1}{2}}^{n+1}-\G_{i,j-\frac{1}{2}}^{n+1}\right).
\end{gather*}
Performing a discrete integration by parts and using zero flux boundary conditions, we have
\begin{align*}
\frac{\Delta t}{\varepsilon}\sum_{j\in\J}v_{j}\left(\G_{i,j+\frac{1}{2}}^{n+1}-\G_{i,j-\frac{1}{2}}^{n+1}\right)&=-\frac{\Delta t}{\varepsilon}\sum_{j\in\J^*}(v_{j+1}-v_{j})\G_{i,j+\frac{1}{2}}^{n+1}\\
&=-\frac{\Delta t}{\varepsilon}\,\Delta x_{i}\sum_{j\in\J^*\setminus\{-L,L\}}\M_{j+\frac{1}{2}}^*\,\left(\frac{f_{i,j+1}^{n+1}}{\M_{j+1}}-\frac{f_{ij}^{n+1}}{\M_{j}}\right).
\end{align*}
Now, performing again a discrete integration by parts and using that the discrete Maxwellian vanishes at endpoints of the velocity domain (see assumption \eqref{hyp:MaxwFP}), we obtain
\begin{align*}
\frac{\Delta t}{\varepsilon}\sum_{j\in\J}v_{j}\left(\G_{i,j+\frac{1}{2}}^{n+1}-\G_{i,j-\frac{1}{2}}^{n+1}\right)&=\frac{\Delta t}{\varepsilon}\,\Delta x_{i}\sum_{j\in\J}(\M_{j+\frac{1}{2}}^*-\M_{j-\frac{1}{2}}^*)\frac{f_{ij}^{n+1}}{\M_{j}}\\
&= -\frac{\Delta t}{\varepsilon}\,\Delta x_{i}\sum_{j\in\J}\Delta v_{j}v_{j}f_{ij}^{n+1},
\end{align*}
which finally yields the result using definition \eqref{eq:def_jS} of $J_{i}^{n+1}$.

In the BGK case, we also multiply \eqref{eq:scheme_f_BGK} by $v_{j}$ and sum over $j$ to obtain
\begin{gather*}
\varepsilon^2\Delta x_{i}(J_{i}^{n+1}-J_{i}^{n})+\Delta t(S_{i+\frac{1}{2}}^{n+1}-S_{i-\frac{1}{2}}^{n+1})+\Delta tm_{2}^{\Delta v}(\rho_{i+\frac{1}{2}}^{n+1}-\rho_{i-\frac{1}{2}}^{n+1})=\\
\frac{\Delta t}{\varepsilon}\Delta x_{i}\sum_{j\in\J}\Delta v_{j}(\rho_{i}^{n+1}\M_{j}-f_{ij}^{n+1})v_{j}.
\end{gather*}
Using \eqref{zero_mean_M}, the right hand side can be written as $-\Delta t\Delta x_{i}J_{i}^{n+1}$, which concludes the proof.
\end{proof}

\begin{thrm}(Asymptotic-preserving property)
 Let us consider a discrete Maxwellian satisfying \eqref{hyp:MaxwFP} (\text{resp.} \eqref{hyp:MaxwBGK}) and let $f_\ve^n = (f_{ij}^n)_{i\in\I,j\in\J}$ for $n\in\NN$ be \Mar{the} solution of  \eqref{eq:scheme_f_FP} (\text{resp.} \eqref{eq:scheme_f_BGK}) with initial data $(f^0_{ij})_{i\in\I,j\in\J}$. Then there is $\rho^n = (\rho_i^n)_{i\in\I}$ for all $n\geq0$ such that when  $\ve\rightarrow0$ one has
 \[
  f_\ve^n \longrightarrow \rho^n \M\quad\text{in}\quad\RR^{\I\times\J}\,, \quad\text{for all}\quad n\geq1
 \]
and the limit satisfies the following finite difference scheme for the heat equation
\begin{equation}
\Delta x_i\frac{\rho^{n+1}_i - \rho^{n}_i}{\Delta t}\ =\ \frac{m_2^{\Delta v}}{2}\left((D_x\rho^{n+1})_{i+1} - (D_x\rho^{n+1})_{i-1}\right)\,,\quad\forall i\in\I
\label{eq:heat_discrete}
\end{equation}
with initial data $\rho_i^0 = \sum_{j\in\J}\Delta v_j f_{ij}^0$.
\label{th:AP}
\end{thrm}

\begin{proof}
Thanks to estimate \eqref{eq:estimmoments_discrete} on the discrete macroscopic density, there is a subsequence $(\ve_k)_k$ such that $\rho_{\ve_k}^n\rightarrow\rho^n$ for all $n\geq 0$.  Thanks to estimate \eqref{eq:L2gamma_discrete} the sequence $(f_{\ve_k}^n)_k\in(\RR^{\I\times\J})^\NN$ converges to $(\rho^n_{i}\,\M_j)_{i\in\I, j\in\J}$ for each $n\geq1$. To identify the limit scheme, we combine \eqref{scheme_rho} and \eqref{scheme_J} to get that for $\ve>0$ one has for all $i\in\I$ and $n\in\NN$ that
\begin{multline*}
\Delta x_i\,(\rho^{n+1}_i - \rho^{n}_i)\ =\ m_2^{\Delta v}\,\Delta t\,\frac{1}{2}\left((D_x\rho^{n+1})_{i+1} - (D_x\rho^{n+1})_{i-1}\right)\\[.5em]
+ \,\Delta t\,\frac{1}{2}\left((D_xS^{n+1})_{i+1} - (D_xS^{n+1})_{i-1}\right) + \ve^2\,(J_{i+\frac{1}{2}}^{n+1}-J_{i-\frac{1}{2}}^{n+1}) - \ve^2\,(J_{i+\frac{1}{2}}^{n}-J_{i-\frac{1}{2}}^{n})\,.
\end{multline*}
Then thanks to estimate \eqref{eq:estimmoments_discrete} it is straightforward to show that the terms of the second line converge to $0$. Finally observe that the limit scheme has a unique solution. Indeed, it consists in inverting a matrix that is strictly diagonally dominant and thus \Mar{that} is non singular. Therefore the whole sequence \Mar{converges}.  
\end{proof}

\begin{rem}
Observe that the scheme for the heat equation keeps a trace of the velocity discretization in the diffusion coefficient $m_2^{\Delta v}$.
\end{rem}

\subsection{Numerical hypocoercivity}

We shall now adapt to the discrete framework the hypocoercivity method proposed by Dolbeault, Mouhot and Schmeiser in \cite{Dolbeault2015} and summarized in Section \ref{sec:hypocoercivity}. The main result of the section is the following.

\begin{thrm}
 Let us consider a discrete Maxwellian satisfying \eqref{hyp:MaxwFP} (\text{resp.} \eqref{hyp:MaxwBGK}) and assume that the number of points $N$ in the space discretization is odd. Then there are constants $C\geq 1$ and $\beta>0$ such that for all $\ve\in(0,1)$, all $\Delta t\leq \Delta t_{\max}$ and all initial data $(f^0_{ij})_{i\in\I,j\in\J}$,  the solution $(f^n_{ij})_{i\in\I,j\in\J,n\in\NN}$ of \eqref{eq:scheme_f_FP}  (\text{resp.} \eqref{eq:scheme_f_BGK}) satisfies
 \begin{equation}
 \|f^n-\mu_f\M\|_{2,\gamma}\ \leq\ C\,\|f^0-\mu_f\M\|_{2,\gamma}\,e^{-\frac{\beta}{2}\,t^n}\,.
  \label{eq:hypocoercivity_discrete}
 \end{equation}
 Moreover, the constants $C$ and $\beta$ do not depend on the size of the discretization $\delta$, and $\Delta t_{\max}>0$ can be chosen arbitrarily.
 \label{th:hypocoercivity_discrete}
\end{thrm}

Once again we prove the theorem assuming that the initial data has zero mean
 \[
  \sum_{(i,j)\in\I\times\J}\Delta x_i\,\Delta v_j\, f^0_{ij}\ =\ 0\,.
 \]
 In the general case one can just replace $f^0_{ij}$ by $f^0_{ij}-\mu_f\M_j$, where $\mu_f=\sum_{(i,j)\in\I\times\J}\Delta x_i\Delta v_j f_{ij}^0$, and correspondingly $f^n_{ij}$ and $\rho^n_i$ by respectively $f^n_{ij} -\mu_f\M_j$ and $\rho^n_i -\mu_f$ in the following arguments.

For all $n\geq 1$, let us define the \emph{discrete modified entropy functional} which reads
\begin{equation}\label{def_H_dis}
\HHH(f^n)\ :=\ \frac{1}{2}\|f^n\|_{2,\gamma}^2+\eta\varepsilon^2\sum_{i\in\I}\Delta x_{i}J_{i}^{n}(D_{x}\phi)_{i}^{n}+\frac{\eta\varepsilon^2}{2}\sum_{i\in\I}\Delta x_{i}\frac{\left((D_{x}\phi)_{i}^{n}-(D_{x}\phi)_{i}^{n-1}\right)^2}{\Delta t},
\end{equation}
where $\eta>0$ will be determined later and $(\phi_{i}^{n})_{i\in\I}$ is the solution of the following discrete Poisson equation
\begin{align}
&-\frac{(D_{x}\phi)_{i+1}^{n}-(D_{x}\phi)_{i-1}^{n}}{2}\ =\ \Delta x_{i}\rho_{i}^n,\quad\forall i\in\I,\label{scheme_phi}\\
&\sum_{i\in\I}\Delta x_{i}\phi_i^n\ =\ 0.\label{cond_phi_mean0}
\end{align}
\begin{rem}
 Observe that the discrete version of the modified entropy has an additional third term compared to \eqref{eq:modifiedEntropy}. It is of order $O(\Delta t)$ and thus consistent with $0$, so that taking limits in the discretization parameters, we recover (at least formally at this stage) the continuous modified entropy \eqref{eq:modifiedEntropy}. It does not perturb the adaptation of the strategy of Section~\ref{sec:hypocoercivity}.
\end{rem}

Let us first justify existence and unicity of $(\phi_{i})_{i\in\I}$ satisfying \eqref{scheme_phi}-\eqref{cond_phi_mean0}. System \eqref{scheme_phi} can be written under the matrix form $A\phi^{n}=\rho^n$, where the rank of the $N\times N$ matrix $A$ is either $N-1$ if $N$ is odd, or $N-2$ if $N$ is even. Observe that $A$ is symmetric with respect to the inner product $\lla\phi,\psi\rra = \sum_{i\in\I}\phi_i\,\psi_i\,\Delta x_i$. In the odd case, the condition \eqref{cond_phi_mean0} corresponds exactly to $(\phi_{i})_{i\in\I}$ being orthogonal (for the latter inner product) to the kernel of $A$ spanned by $(1,\dots,1)\in\RR^\I$. Hence provided that $N$ is odd, which is assumed in Theorem~\ref{th:hypocoercivity_discrete}, the system \eqref{scheme_phi} supplemented with zero mean value condition \eqref{cond_phi_mean0} admits a unique solution.

Let us now derive the discrete estimates on $(\phi_{i})_{i\in\I}$.

\begin{lem}
 Under the assumptions of Theorem~\ref{th:hypocoercivity_discrete}, one has for all $n\in\NN$
 \begin{align}
  \|D_x\phi^n\|_{2}&\ \leq\ C_P\,\|\rho^n\|_{2}\,,\label{eq:dxphi_discrete}\\[1em]
  \left\|\frac{(D_{x}\phi)_{i}^{n+1}-(D_{x}\phi)_{i}^{n}}{\Delta t}\right\|_{2}&\ \leq\ \|J_i^{n+1}\|_{2}\,,\label{eq:dtdxphi_discrete}
 \end{align}
 where $C_P$ is the discrete Poincaré constant of Lemma~\ref{lem:poincare_discrete}.
\end{lem}
\begin{proof}
Multiplying the discrete Poisson equation \eqref{scheme_phi} by $\phi_{i}^n$, summing on $i\in\I$ and performing a discrete integration by parts, we obtain
\[
-\sum_{i\in\I}\Delta x_{i}\left((D_{x}\phi)_{i}^{n}\right)^2=-\sum_{i\in\I}\Delta x_{i}\rho_{i}^{n}\phi_{i}^{n},
\]
which yields thanks to Cauchy-Schwarz inequality
\begin{equation*}
\sum_{i\in\I}\Delta x_{i}\left((D_{x}\phi)_{i}^{n}\right)^2\leq \left(\sum_{i\in\I}\Delta x_{i}(\rho_{i}^n)^2\right)^{\frac{1}{2}}\left(\sum_{i\in\I}\Delta x_{i}(\phi_{i}^n)^2\right)^{\frac{1}{2}}.
\end{equation*}
Using the Poincaré inequality \eqref{ineq_Poincare_x} we get the first estimate.

Now, considering scheme \eqref{scheme_phi} at time $n$ and $n+1$, we can write
\begin{equation*}
-\frac{1}{2}\left[\left((D_{x}\phi)_{i+1}^{n+1}-(D_{x}\phi)_{i-1}^{n+1}\right)-\left((D_{x}\phi)_{i+1}^{n}-(D_{x}\phi)_{i-1}^{n}\right)\right]\ =\ \Delta x_{i}(\rho_{i}^{n+1}-\rho_{i}^n).
\end{equation*}
Multiplying this equality by $\phi_{i}^{n+1}-\phi_{i}^n$, summing on $i\in\I$ and integrating by parts, we get
\begin{equation*}
\sum_{i\in\I}\Delta x_{i}\left[(D_{x}\phi)_{i}^{n+1}-(D_{x}\phi)_{i}^{n}\right]^2=\sum_{i\in\I}\Delta x_{i}(\rho_{i}^{n+1}-\rho_{i}^n)(\phi_{i}^{n+1}-\phi_{i}^n).
\end{equation*}
Using the scheme \eqref{scheme_rho}, performing a discrete integration by parts and applying Cauchy-Schwarz inequality, we obtain:
\begin{align*}
\sum_{i\in\I}\Delta x_{i}\left[(D_{x}\phi)_{i}^{n+1}-(D_{x}\phi)_{i}^{n}\right]^2&=-\Delta t\sum_{i\in\I}\frac{J_{i+1}^{n+1}-J_{i-1}^{n+1}}{2}(\phi_{i}^{n+1}-\phi_{i}^n)\\
&=\Delta t\sum_{i\in\I}\Delta x_{i} J_{i}^{n+1}\left[(D_{x}\phi)_{i}^{n+1}-(D_{x}\phi)_{i}^{n}\right]\\
&\leq \Delta t\left(\sum_{i\in\I}\Delta x_{i}(J_{i}^{n+1})^2\right)^{\frac{1}{2}}\left(\sum_{i\in\I}\Delta x_{i}\left[(D_{x}\phi)_{i}^{n+1}-(D_{x}\phi)_{i}^{n}\right]^2\right)^{\frac{1}{2}},
\end{align*}
from which we deduce \eqref{eq:dtdxphi_discrete}.
\end{proof}

Now we show that for $\eta$ small enough the modified entropy functional \eqref{def_H_dis} is an equivalent $\|\cdot\|_{2,\gamma}$ norm uniformly for small $\ve$ and mesh size. 

\begin{lem}
Under the assumptions of Theorem~\ref{th:hypocoercivity_discrete}, and assuming that $\varepsilon\leq 1$ and $\Delta t\leq \Delta t_{max}$, one has for all $n\geq 1$
\begin{equation}\label{eq:compare_H_E}
\left(\frac{1}{2}-\eta\,(m_2^{\Delta v})^{1/2}\right)\|f^n\|_{2,\gamma}^2\ \leq\ \HHH(f^n)\ \leq\ \left(\frac{1}{2}+\eta\,(m_2^{\Delta v})^{1/2}C_P+\frac{\eta}{2}m_2^{\Delta v}\Delta t_{max}\right)\|f^n\|_{2,\gamma}^2.
\end{equation}
\end{lem}

\begin{proof}
On the one hand, using the Cauchy-Schwarz inequality, we have
\[
\eta\varepsilon^2\left|\sum_{i\in\I}\Delta x_{i}J_{i}^n(D_{x}\phi)_{i}^n\right|\leq \eta\varepsilon^2\left(\sum_{i\in\I}\Delta x_{i}(J_{i}^n)^2\right)^{\frac{1}{2}}\left(\sum_{i\in\I}\Delta x_{i}\left((D_{x}\phi)_{i}^n\right)^2\right)^{\frac{1}{2}}\,,
\]
so that with \eqref{eq:estmomentsLinfj_discrete}, \eqref{eq:dxphi_discrete} and \eqref{eq:estmomentsLinfrho_discrete} one has
\[
 \eta\varepsilon^2\left|\sum_{i\in\I}\Delta x_{i}J_{i}^n(D_{x}\phi)_{i}^n\right|\ \leq\ \eta\,(m_2^{\Delta v})^{1/2}\,C_P\,\ve\,\|f^n\|_{2,\gamma}^2.
\]
On the other hand the last term of $\HHH(f^n)$ may be estimated using \eqref{eq:dtdxphi_discrete} and \eqref{eq:estmomentsLinfj_discrete}% and \eqref{eq:estimmoments_discrete}, we have:
\[
\frac{\eta\varepsilon^2}{2}\,\sum_{i\in\I}\Delta x_{i}\frac{\left((D_{x}\phi)_{i}^n-(D_{x}\phi)_{i}^{n-1}\right)^2}{\Delta t}\ \leq\ \frac{\eta}{2}\,m_2^{\Delta v}\,\Delta t\,\|f^{n+1}\|_{2,\gamma}^2 .
\]
We finally establish \eqref{eq:compare_H_E} using that $\varepsilon\leq 1$ and $\Delta t\leq \Delta t_{max}$.
\end{proof}

\begin{prop}
Under the assumptions of Theorem~\ref{th:hypocoercivity_discrete}, there is $\eta_2>0$ such that for all $\varepsilon\leq 1$, $\Delta t\leq \Delta t_{max}$ and $\eta\leq\eta_2$, 
\begin{equation*}
\frac{\HHH(f^{n+1})-\HHH(f^n)}{\Delta t}\,+\, K(\eta)\left(\|f^{n+1} - \rho^{n+1}\M\|_{2,\gamma}^2+\|\rho^{n+1}\|_2^2\right)\ \leq\ 0\,,\qquad\forall n\geq 1\,\,,
\end{equation*}
with $K(\eta) = \frac{1}{2}\min\left(1-\eta \overline{m}_{2},\eta \underline{m}_{2}\right)$.
\end{prop}

 \begin{proof}
 Let us estimate the discrete time derivative of $\HHH$. For all $n\geq 1$, we have
\begin{equation}\label{estimHD1}
\HHH(f^{n+1})-\HHH(f^{n})\ =\ \frac{1}{2}(\|f^{n+1}\|_{2,\gamma}^2-\|f^{n}\|_{2,\gamma}^2)+\eta\varepsilon^2 T_{1}^{n}+\eta\varepsilon^2 T_{2}^{n},
\end{equation}
with
\begin{align*}
&T_{1}^{n}=\sum_{i\in\I}\Delta x_{i}\left(J_{i}^{n+1}(D_{x}\phi)_{i}^{n+1}-J_{i}^{n}(D_{x}\phi)_{i}^{n}\right),\\
&T_{2}^n=\frac{1}{2}\sum_{i\in\I}\frac{\Delta x_{i}}{\Delta t}\left(\left((D_{x}\phi)_{i}^{n+1}-(D_{x}\phi)_{i}^{n}\right)^2-\left((D_{x}\phi)_{i}^{n}-(D_{x}\phi)_{i}^{n-1}\right)^2\right).
\end{align*}
 Concerning the first term, we proved in Lemma~\ref{lem:basic_estimates_discrete} that
\begin{equation}\label{eq:estimEI_disbis}
\frac{1}{2}(\|f^{n+1}\|_{2,\gamma}^2-\|f^{n}\|_{2,\gamma}^2)\leq -\frac{\Delta t}{\varepsilon^2}\|f^{n+1} - \rho^{n+1}\M\|_{2,\gamma}^2.
\end{equation}
Then, we can write:
\begin{equation}\label{eq:T1_2}
T_{1}^{n}=\sum_{i\in\I}\Delta x_{i}(J_{i}^{n+1}-J_{i}^n)(D_{x}\phi)_{i}^{n+1}+\sum_{i\in\I}\Delta x_{i}J_{i}^{n+1}\left((D_{x}\phi)_{i}^{n+1}-(D_{x}\phi)_{i}^{n}\right)+R^{n},
\end{equation}
where 
\begin{equation*}
R^n:=\sum_{i\in\I}\Delta x_{i}\left(J_{i}^{n}(D_{x}\phi)_{i}^{n+1}-J_{i}^{n+1}(D_{x}\phi)_{i}^{n+1}+J_{i}^{n+1}(D_{x}\phi)_{i}^{n}-J_{i}^n(D_{x}\phi)_{i}^{n}\right).
\end{equation*}
We remark that the two first sums in the right hand side of \eqref{eq:T1_2} correspond to a discrete \Mar{version} of $(\lla\partial_tj^\ve,\partial_x\phi^\ve\rra+\lla j^\ve,\partial_t\partial_x\phi^\ve\rra)\Delta t$ .

Let us now study the remainder term $R^{n}$. We have
\begin{equation*}
R^n=\sum_{i\in\I}\Delta x_{i}(J_{i}^n-J_{i}^{n+1})\left((D_{x}\phi)_{i}^{n+1}-(D_{x}\phi)_{i}^{n}\right),
\end{equation*}
which after a discrete integration by parts leads to
\begin{equation*}
R^n=\sum_{i\in\I}\Delta x_{i}\left((D_{x}J)_{i}^{n+1}-(D_{x}J)_{i}^{n}\right)(\phi_{i}^{n+1}-\phi_{i}^{n}).
\end{equation*}
Using scheme \eqref{scheme_rho}, we get
\begin{equation*}
R^n=-\frac{1}{\Delta t}\sum_{i\in\I}\Delta x_{i}(\rho_{i}^{n+1}-2\rho_{i}^n+\rho_{i}^{n-1})(\phi_{i}^{n+1}-\phi_{i}^{n}),
\end{equation*}
and then using scheme \eqref{scheme_phi} followed by a discrete integration by parts, we obtain:
\begin{align*}
R^{n}&=\frac{1}{2\Delta t}\sum_{i\in\I}\left[\left((D_{x}\phi)_{i+1}^{n+1}-(D_{x}\phi)_{i-1}^{n+1}\right)-2\left((D_{x}\phi)_{i+1}^{n}-(D_{x}\phi)_{i-1}^{n}\right)\right.\\
&\quad\quad\quad\left.+\left((D_{x}\phi)_{i+1}^{n-1}-(D_{x}\phi)_{i-1}^{n-1}\right)\right](\phi_{i}^{n+1}-\phi_{i}^{n})\\
&=-\frac{1}{\Delta t}\sum_{i\in\I}\Delta x_{i}\left((D_{x}\phi)_{i}^{n+1}-2(D_{x}\phi)_{i}^{n}+(D_{x}\phi)_{i}^{n-1}\right)\left((D_{x}\phi)_{i}^{n+1}-(D_{x}\phi)_{i}^{n}\right).
\end{align*}
Using that $-a(a-b) + (a^2-b^2)/2 = -(a-b)^2/2\leq 0$ with  $a=(D_{x}\phi)_{i}^{n+1}-(D_{x}\phi)_{i}^{n}$ and $b=(D_{x}\phi)_{i}^{n}-(D_{x}\phi)_{i}^{n-1}$, we finally get
\begin{equation}\label{eq:estimR}
R^n + T_{2}^n\leq 0\,.
\end{equation}
Then gathering \eqref{eq:estimEI_disbis} and \eqref{eq:estimR} in \eqref{estimHD1} results in
\begin{equation}\label{estimHD2}
\HHH(f^{n+1})-\HHH(f^{n})\leq -\frac{\Delta t}{\varepsilon^2}\|f^{n+1} - \rho^{n+1}\M\|_{2,\gamma}^2+T_{3}^n,
\end{equation}
with
\begin{equation*}
T_{3}^n:=\eta\varepsilon^2\sum_{i\in\I}\Delta x_{i}(J_{i}^{n+1}-J_{i}^n)(D_{x}\phi)_{i}^{n+1}+\eta\varepsilon^2\sum_{i\in\I}\Delta x_{i}J_{i}^{n+1}\left((D_{x}\phi)_{i}^{n+1}-(D_{x}\phi)_{i}^{n}\right).
\end{equation*}

It remains now to estimate this term $T_{3}^n$. We proceed as in the continuous framework. Using scheme \eqref{scheme_J} in the first sum of the right hand side, we write:
\begin{equation*}
T_{3}^n=T_{31}^{n}+T_{32}^{n}+T_{33}^n+T_{34}^n,
\end{equation*}
where
\begin{align*}
T_{31}^n&:=-\eta\Delta t\sum_{i\in\I}\Delta x_{i}(D_{x}S)_{i}^{n+1}(D_{x}\phi)_{i}^{n+1},\\
T_{32}^{n}&:=-\eta\Delta t\sum_{i\in\I}\Delta x_{i}m_{2}^{\Delta v}(D_{x}\rho)_{i}^{n+1}(D_{x}\phi)_{i}^{n+1},\\
T_{33}^{n}&:=-\eta\Delta t\sum_{i\in\I}\Delta x_{i}J_{i}^{n+1}(D_{x}\phi)_{i}^{n+1},\\
T_{34}^n&:=\eta\varepsilon^2\sum_{i\in\I}\Delta x_{i}J_{i}^{n+1}\left((D_{x}\phi)_{i}^{n+1}-(D_{x}\phi)_{i}^{n}\right).
\end{align*}
Applying a discrete integration by parts and using scheme \eqref{scheme_phi}, we have
\[
T_{31}^{n}=\eta\Delta t\sum_{i\in\I}S_{i}^{n+1}\frac{(D_{x}\phi)_{i+1}^{n+1}-(D_{x}\phi)_{i-1}^{n+1}}{2}=-\eta\Delta t\sum_{i\in\I}\Delta x_{i}S_i^{n+1}\rho_{i}^{n+1}\,.
\]
Then the Cauchy-Schwarz inequality and \eqref{eq:estmomentsL2S_discrete} yields
\begin{equation}\label{estim_T31}
|T_{31}^n|\leq \eta\Delta t\|S^{n+1}\|_2\|\rho^{n+1}\|_2\leq \eta\Delta t\varepsilon\sqrt{\overline{m}_{4}-\underline{m}_{2}^2}\frac{\|f^{n+1} - \rho^{n+1}\M\|_{2,\gamma}}{\varepsilon}\|\rho^{n+1}\|_2\,.
\end{equation}

Performing a discrete integration by parts and using again scheme \eqref{scheme_phi}, we get
\begin{equation}\label{estim_T32}
T_{32}^n\ =\ -\eta\Delta t m_{2}^{\Delta v}\sum_{i\in\I}\Delta x_{i}(\rho_{i}^{n+1})^2\ =\ -\eta\,\Delta t\,m_{2}^{\Delta v}\,\|\rho^{n+1}\|_2^2\,.
\end{equation}

Then, applying Cauchy-Schwarz inequality followed by \eqref{eq:estmomentsL2j_discrete} and \eqref{eq:dxphi_discrete} we get
\begin{equation}\label{estim_T33}
 |T_{33}^n|\leq \eta\,\Delta t\,\|J^{n+1}\|_2\,\|D_x\phi^{n+1}\|_2\leq \eta\,\Delta t\,C_P\,\sqrt{\overline{m}_{2}}\,\frac{\|f^{n+1} - \rho^{n+1}\M\|_{2,\gamma}}{\varepsilon}\|\rho^{n+1}\|_2\,.
\end{equation}

The last remainder term is estimated again using Cauchy-Schwarz inequality followed by \eqref{eq:estmomentsL2j_discrete} and \eqref{eq:dtdxphi_discrete}, yielding
\begin{equation}\label{estim_T34}
|T_{34}^{n}|\leq \eta\varepsilon^2\,\|J^{n+1}\|_2\,\|D_x\phi^{n+1} - D_x\phi^{n}\|_2\,\leq \eta\Delta t \overline{m}_{2}\|f^{n+1} - \rho^{n+1}\M\|_{2,\gamma}^2.
\end{equation}

Gathering \eqref{estim_T31}, \eqref{estim_T32}, \eqref{estim_T33} and \eqref{estim_T34} in \eqref{estimHD2} then gives:
\begin{gather*}
\HHH(f^{n+1})-\HHH(f^{n})+\Delta t (\frac{1}{\varepsilon^2}-\eta \overline{m}_{2})\|f^{n+1} - \rho^{n+1}\M\|_{2,\gamma}^2+\Delta t\eta \underline{m}_{2}\|\rho^{n+1}\|_2^2\\
\leq \Delta t\eta( \sqrt{\overline{m}_{4}-\underline{m}_{2}^2}+\frac{C_P\,\sqrt{\overline{m}_{2}}}{\varepsilon})\|f^{n+1} - \rho^{n+1}\M\|_{2,\gamma}\|\rho^{n+1}\|_2.\nonumber
\end{gather*}

It follows that for any $\eta\in (0,\eta_2)$ with 
\begin{equation*}
\eta_2:=\frac{\underline{m}_{2}}{(\sqrt{\overline{m}_{4}-\underline{m}_{2}^2}+C_P\,\sqrt{\overline{m}_{2}})^2+\underline{m}_{2}\overline{m}_{2}},
\end{equation*}
and $\varepsilon\in (0,1)$, one has
\begin{equation*}
    \HHH(f^{n+1})-\HHH(f^n)+\Delta t K(\eta)\left( \|f^{n+1}-\rho^{n+1}\MM\|_{2,\gamma}^2+\|\rho^{n+1}\|_{2}^2  \right)\leq 0,
\end{equation*}
with 
\begin{equation*}
K(\eta)=\frac{1}{2}\min\left(1-\eta \overline{m}_{2},\eta \underline{m}_{2}\right).
\end{equation*}

\end{proof}

\begin{proof}[Proof of Theorem~\ref{th:hypocoercivity_discrete}]
To conclude the proof \Max{we observe that $f^{n+1}- \rho^{n+1}\M$ and $\rho^{n+1}\M$ are orthogonal with respect to $\|\cdot\|_{2,\gamma}$} and thus with \eqref{eq:compare_H_E} one obtains
%we use the triangle inequality, \Mar{$(a+b)^2\leq 2(a^2 + b^2)$} followed by  to obtain
\begin{equation*}
\HHH(f^{n+1})-\HHH(f^n)+\Delta t \kappa(\eta)\,\HHH(f^{n+1})\leq 0\quad\forall n\geq 1,
\end{equation*}
and $\kappa(\eta) = \Mar{2}K(\eta)/K_2(\eta)$ where $K_2(\eta) = 1+2\eta\sqrt{\overline{m}_2}C_P+\eta\overline{m}_2\Delta t_{\max}$. It implies 
\[
\HHH(f^n)\ \leq\ \HHH(f^{1})\,\left(1+\Delta t\,\kappa(\eta)\right)^{-(n-1)}\, = \HHH(f^1) \exp(-t^{n-1}(\Delta t)^{-1}\log\left(1+\Delta t\,\kappa(\eta)\right))\,.
\]
Since $s\mapsto \log\left(1+s \kappa(\eta)\right)/s$ is a non increasing function on $]0,+\infty[$ and $\Delta t\leq \Delta t_{max}$, one has 
\[\HHH(f^{n+1})\leq e^{\beta \Delta t_{max}}\HHH(f^{1})e^{-\beta t^{n}} \quad \forall n\geq 1,\]
with maximal uniform rate
\begin{equation*}
 \beta = \frac{1}{\Delta t_{max}}\log \left(1+\Delta t_{max}\kappa(\eta)\right).
\end{equation*}
Finally \[\HHH(f^{1})\leq \frac{K_{2}(\eta)}{2}\|f^1\|^2_{2,\gamma}\leq \frac{K_{2}(\eta)}{2}\|f^0\|^2_{2,\gamma}\,.\]

Provided that $0<\eta<\min(\eta_1,\eta_2)$, with $\eta_1 = 1/(2\underline{m}_2^{1/2})$ to ensure that the constant at the  left hand side of \eqref{eq:compare_H_E} is positive, one finally obtains \eqref{eq:hypocoercivity_discrete}, which concludes the proof of Theorem~\ref{th:hypocoercivity_discrete}.

 \end{proof}

\section{Implementation}\label{sec:implementation}

In practice, the implementation of the schemes written as in \eqref{eq:scheme_f_FP} and \eqref{eq:scheme_f_BGK} leads to two numerical issues. First, the matrix that should be inverted is ill-conditioned for small $\ve$. This leads to large numerical errors near the diffusive regime and deteriorates the AP property of the schemes. Second, the evaluation of integrated differences that are of the same order, such as $\|f^n-\mu_f\M\|_{2,\M}$ when $n$ is large leads to the accumulations of machine epsilons and poor approximations of these quantities.

In this section, we present an equivalent formulation of the schemes \eqref{eq:scheme_f_FP} and \eqref{eq:scheme_f_BGK} which solves both these issues. It is based on a perturbative micro-macro  formulation introduced in Section~\ref{sec:micmac}. In the next subsection, we justify the uniform conditioning of the reformulated linear systems on the whole range  of $\ve$ (including $\ve =0$). In Section~\ref{sec:matrix}, we write the linear systems in explicit matrix form so that the reader may implement the scheme easily and replicate the numerical results of Section~\ref{sec:numres}. The test cases illustrate the theoretical results of Theorem~\ref{th:AP} and Theorem~\ref{th:hypocoercivity_discrete}, namely the asymptotic accuracy of our schemes in both the diffusion limit and the large time asymptotics. 

\subsection{Micro-macro formulations of the schemes}\label{sec:micmac}

Given $(f_{ij}^n)_{i\in\I, j\in\J, n\in\NN}$, let us  introduce the micro and macro unknowns, $(h_{ij}^n)_{i\in\I, j\in\J, n\in\NN}$ and $(\lambda_i^n)_{i\in\I, n\in\NN}$ respectively, which satisfy for all $i\in\I$, $j\in\J$ and $n\in\NN$
 \begin{align}
 f_{ij}^n&= \mu_f\,\M_j  + \lambda_i^n\,\M_j + \ve\,h_{ij}^n\,\M_j\,,\label{eq:micromacro}\\[.5em]
 \lambda_{i}^n&= \rho_i^n - \mu_f = \sum_{j\in\I} f_{ij}^n\,\Delta v_j-\mu_f\label{eq:lambda}\,.
 \end{align}
An immediate consequence of these expressions is that 
\begin{equation}
\sum_{j\in\J}h_{ij}^n\,\M_j\,\Delta v_j\ =\ 0\,\qquad\forall i\in\I\,,\forall n\in\NN\,.
 \label{eq:h_meanfree}
\end{equation}
If one writes the counterpart of the continuity equation \eqref{scheme_rho} with the new unknowns, one gets the following macroscopic evolution equation
\begin{equation}
 \Delta x_i(\lambda_i^{n+1} - \lambda_i^{n})\,+\,\Delta t\,\sum_{j\in\J}\,\Delta v_j\,v_j\,\M_j\,\frac{h_{i+1,j}^{n+1} - h_{i-1,j}^{n+1}}{2}\ =\ 0\,,\qquad\forall i\in\I\,.
 \label{eq:scheme_lambda}
\end{equation}
\begin{itemize}

\item \emph{In the Fokker-Planck case}, if one plugs \eqref{eq:micromacro} into the scheme \eqref{eq:scheme_f_FP} and then \Mar{subtracts} \eqref{eq:scheme_lambda} multiplied by $-\Delta v_j\,\ve$, one obtains the following microscopic evolution equation
\begin{multline}
 \,\ve^2\,\Delta x_i\,\Delta v_j\,(h_{ij}^{n+1} - h_{ij}^{n}) \,+\, \Delta t\,\Delta v_j\,v_j\frac{\lambda_{i+1}^{n+1}-\lambda_{i-1}^{n+1}}{2}\\
 + \ve\,\Delta t\,\Delta v_j\,\left(v_j\,\frac{h_{i+1,j}^{n+1} - h_{i-1,j}^{n+1}}{2} - \sum_{k\in\J}\,\Delta v_k\,v_k\,\M_k\,\frac{h_{i+1,k}^{n+1} - h_{i-1,k}^{n+1}}{2}\right)\\
 \ =\ \Delta t\,\Delta x_i\left(\frac{\M_{j+1/2}^*}{\Delta v_{j+1/2}\,\M_j}\left(h_{i,j+1}^{n+1} - h_{ij}^{n+1}\right) - \frac{\M_{j-1/2}^*}{\Delta v_{j-1/2}\,\M_j}\left(h_{ij}^{n+1} - h_{i,j-1}^{n+1}\right)\right)\,.
 \label{eq:scheme_h_FP}
\end{multline}

\item \emph{In the BGK case}, if one plugs \eqref{eq:micromacro} into the scheme \eqref{eq:scheme_f_BGK}, one obtains the following microscopic evolution equation
\begin{multline}
 \,\ve^2\,\Delta x_i\,\Delta v_j\,(h_{ij}^{n+1} - h_{ij}^{n}) \,+\, \Delta t\,\Delta v_j\,v_j\frac{\lambda_{i+1}^{n+1}-\lambda_{i-1}^{n+1}}{2}\\
 + \ve\,\Delta t\,\Delta v_j\,\left(v_j\,\frac{h_{i+1,j}^{n+1} - h_{i-1,j}^{n+1}}{2} - \sum_{k\in\J}\,\Delta v_k\,v_k\,\M_k\,\frac{h_{i+1,k}^{n+1} - h_{i-1,k}^{n+1}}{2}\right)\\
 \ =\ -\Delta t\,\Delta x_i\,\Delta v_j  h_{ij}^{n+1}\,.
 \label{eq:scheme_h_BGK}
\end{multline}
\end{itemize}

Based on the previous computations, one can define three schemes.  In all three cases we start from the given initial data $(f_{ij}^0)_{i, j}\in\RR^{\I\times\J}$ and define  $(h_{ij}^0)_{i, j}\in\RR^{\I\times\J}$ and \Mar{$(\lambda_i^0)_{i }\in\RR^{\I}$} according to \eqref{eq:micromacro} and \eqref{eq:lambda}. If $\ve=0$, we define $h_{ij}^0 = 0$ for all $i,j$. Then the schemes are as follows.
\begin{itemize}
 \item[$(\mathbf{S}_{f}^\ve)$] The unknown $(f_{ij}^n)_{i, j, n}\in\RR^{\I\times\J\times\NN}$ solves \eqref{eq:scheme_f_FP} in the Fokker-Planck case and \eqref{eq:scheme_f_BGK} in the BGK case. The micro and macro unknowns $(h_{ij}^n)_{i, j, n}\in\RR^{\I\times\J\times\NN}$ and $(\lambda_i^n)_{i, n}\in\RR^{\I\times\NN}$ are defined by \eqref{eq:micromacro} and \eqref{eq:lambda}.
 \item[$(\mathbf{S}_{\lambda,h}^\ve)$] The micro and macro unknowns $(h_{ij}^n)_{i, j, n}$ and $(\lambda_i^n)_{i, n}$ solve \eqref{eq:scheme_lambda}, \eqref{eq:scheme_h_FP} in the Fokker-Planck case and \eqref{eq:scheme_lambda}, \eqref{eq:scheme_h_BGK} in the BGK case. The unknown $(f_{ij}^n)_{i, j, n}$ is defined by \eqref{eq:micromacro}.
 \item[$(\tilde{\mathbf{S}}_{\lambda,h}^\ve)$] The micro and macro unknowns $(h_{ij}^n)_{i, j, n}$ and $(\lambda_i^n)_{i, n}$ solve the overdetermined linear system \eqref{eq:h_meanfree}, \eqref{eq:scheme_lambda}, \eqref{eq:scheme_h_FP} in the Fokker-Planck case and \eqref{eq:h_meanfree}, \eqref{eq:scheme_lambda}, \eqref{eq:scheme_h_BGK} in the BGK case. The unknown $(f_{ij}^n)_{i, j, n}$ is defined by \eqref{eq:micromacro}.
\end{itemize}

These schemes happen to be well-posed and equivalent in the kinetic regime $\ve>0$ but not in the diffusive regime $\ve=0$.

\begin{prop}The three versions of the schemes satisfy the following properties.
 \begin{itemize}
  \item[(i)] If $\ve>0$, $(\mathbf{S}_{f}^\ve)$, $(\mathbf{S}_{\lambda,h}^\ve)$ and $(\tilde{\mathbf{S}}_{\lambda,h}^\ve)$ define the same unique $(f_{ij}^n)_{i, j, n}$, $(h_{ij}^n)_{i, j, n}$ and $(\lambda_i^n)_{i, n}$.
  \item[(ii)] If $\ve=0$, the equivalence between the schemes does not hold. More precisely, $(\mathbf{S}_{f}^0)$ is ill-posed, $(\mathbf{S}_{\lambda,h}^0)$ is well-posed in the BGK case and ill-posed in the Fokker-Planck case and $(\tilde{\mathbf{S}}_{\lambda,h}^0)$ is well-posed. 
 \end{itemize}
\end{prop}
\begin{proof}~\\[-1em]
 \begin{itemize}

  \item[(i)] The computations of the beginning of this section already showed that $(\mathbf{S}_{f}^\ve) \Rightarrow (\mathbf{S}_{\lambda,h}^\ve)$. Let us now prove that $(\mathbf{S}_{\lambda,h}^\ve) \Rightarrow (\tilde{\mathbf{S}}_{\lambda,h}^\ve)$. At $n=0$ one has $\sum_{j\in\J}\Delta v_j\,h_{ij}^{0}\,\M_j\ =\ 0$ for all $i\in\I$. In the Fokker-Planck case,  one can multiply \eqref{eq:scheme_h_FP} by $\M_j$  and \Mar{sum} over $j\in\J$ to obtain 
  \[
   \,\sum_{j\in\J}\Delta v_j\,(\ve^2\,h_{ij}^{n+1} - \ve^2\,h_{ij}^{n})\,\M_j\ =\ 0\,,\qquad \forall i\in\I\,,\quad \forall n\in\NN\,.
  \]
  In the BGK case, it yields 
  \begin{equation}
   \sum_{j\in\J}\Delta v_j\,((\ve^2 + \Delta t)\,h_{ij}^{n+1} - \ve^2\,h_{ij}^{n})\,\M_j\ =\ 0\,,\qquad \forall i\in\I\,,\quad \forall n\in\NN\,.
   \label{eq:meanfree_BGK}
  \end{equation}
  Therefore \eqref{eq:h_meanfree} holds. Finally let us show that $(\tilde{\mathbf{S}}_{\lambda,h}^\ve) \Rightarrow (\mathbf{S}_{f}^\ve)$. 
  Since \eqref{eq:h_meanfree} holds, by multiplying \eqref{eq:micromacro} by $\Delta v_j$  and summing over $j\in\J$ one obtains \eqref{eq:lambda}. Then, just add \eqref{eq:scheme_lambda} multiplied by $\ve\Delta v_j$ to \eqref{eq:scheme_h_FP} (\emph{resp.} \eqref{eq:scheme_h_BGK}) to obtain \eqref{eq:scheme_f_FP} (\emph{resp.} \eqref{eq:scheme_f_BGK}).

  \item[(ii)] Let $\ve = 0$.  In the case of the Fokker-Planck operator, let us define the discrete collision matrix $Q_{FP}^\delta$ satisfying for all $g = (g_j)_{j\in\J}$ 
  \[
   Q_{FP}^\delta\, g\ :=\ \left(\frac{1}{\Delta v_j}\left(\frac{\M_{j+1/2}^*}{\Delta v_{j+1/2}\,\M_j}\left(g_{j+1} - g_{j}\right) - \frac{\M_{j-1/2}^*}{\Delta v_{j-1/2}\,\M_j}\left(g_{j} - g_{j-1}\right)\right)\right)_{j\in\J}\,.
  \]
  Since one has $\lla Q_{FP}^\delta\,g,\,g\rra_{2,\M} = \|D_v g\|_{2,\M}^2$, it is clear that $\mathrm{Ker}(Q_{FP}^\delta) = \mathrm{Span}\{(1,1,\dots,1)\}$ and therefore the linear system \eqref{eq:scheme_h_FP} (with $\ve =0$) has infinitely many solutions since $(v_j)_{j\in\J}\in \mathrm{Ker}(Q_{FP}^\delta)^\perp = \{g,\ \text{s. t.}\ \sum_{j\in\J}g_j\M_j\Delta v_j = 0\}$. If one multiplies  \eqref{eq:scheme_h_FP} by $v_j\M_j$ and sums over $j$, one gets for any solution
  \[
   \sum_{j\in\J}\Delta v_j\,v_j\,\M_j\,h_{ij}^{n+1}\ =\ -m_2^{\Delta v}\,\frac{\lambda_{i+1}^{n+1} - \lambda_{i-1}^{n+1}}{2\Delta x_i}\,.
  \]
  Then by plugging this into \eqref{eq:scheme_lambda}, we infer that $(\lambda_i^{n+1})_{i\in\I}$ is determined uniquely since it satisfies our limit heat equation \eqref{eq:heat_discrete}.
  Hence, the linear system of $(\mathbf{S}_{\lambda,h}^0)$ has infinitely many solutions. If we now consider the scheme $(\tilde{\mathbf{S}}_{\lambda,h}^0)$ one has the additional equations $\sum_{j\in\J}h_{ij}\,\M_j\,\Delta v_j\ =\ 0$ for all $i\in\I$ which correspond to $(h_{ij})_j\in\mathrm{Ker}(Q_{FP}^\delta)^\perp$ for all $i\in\I$. The solution is now unique which proves well-posedness of $(\tilde{\mathbf{S}}_{\lambda,h}^0)$.
  
   In the case of the BGK operator, the situation is simpler since the discrete collision operator $Q_{BGK}^\delta$ is a diagonal positive definite matrix. Hence both $(\mathbf{S}_{\lambda,h}^0)$ and $(\tilde{\mathbf{S}}_{\lambda,h}^0)$ are well-posed. With similar arguments than those developed in the foregoing proof, it is easy to show that $(\mathbf{S}_{f}^0)$ is ill-posed in both the Fokker-Planck and BGK cases.
 \end{itemize}
\end{proof}
An immediate consequence of the previous proposition is the following
\begin{cor}
 The condition number of the linear system in $(\tilde{\mathbf{S}}_{\lambda,h}^\ve)$ is uniformly bounded with respect to $\ve\in[0,1]$.
\end{cor}
\begin{proof}
 By invertibility of the matrices, the condition number is well defined on the whole range $\ve\in[0,1]$. Moreover it is clearly continuous in $\ve$. Hence the boundedness follows from uniform continuity. 
\end{proof}

\begin{rem}
 By taking our schemes for the kinetic equations at $\ve =0$ in their formulations $(\tilde{\mathbf{S}}_{\lambda,h}^0)$, one is exactly solving the discrete heat equation \eqref{eq:heat_discrete} for $(\lambda_i)_{i\in\I}$, as expected in the diffusion limit.
\end{rem}

\subsection{Explicit matrix formulation}\label{sec:matrix}
Following the results of the previous section, our schemes should be implemented in their robust form  $(\tilde{\mathbf{S}}_{\lambda,h}^\ve)$. In this section, we write the explicit matrix formulation of this scheme. At time $t^n$, both schemes rewrite in linear system form \[\mathbb{M}^\varepsilon_\text{col}\,\mathbf{U}^{n+1}\ =\ \mathbb{D}^\varepsilon\,\mathbf{U}^{n}\] with $\text{col}\ =\ FP\text{ or }{BGK}$. The vector $\U^n$ has two sets of components corresponding to the macro and micro variable. They are indexed by control volumes. The first set is $\U^n_{\X_{i}} = \lambda^n_i$ for $i\in\I$ and the second is $\U^n_{K_{ij}} = h^n_{ij}$ with $(i,j)\in\I\times\J$. From the knowledge of $\U^n$, one can compute the usual distribution function $(f^n_{ij})_{i,j}$ thanks to \eqref{eq:micromacro}. Since the linear system is rectangular (but has a unique solution as shown previously), it can be solved thanks to the pseudo-inverse
\begin{equation}
\mathbf{U}^{n+1}\ =\  \left[(\mathbb{M}^\varepsilon_{\text{col}})^\top\,\mathbb{M}^\varepsilon_{\text{col}}\right]^{-1}(\mathbb{M}^\varepsilon_{\text{col}})^\top\,\mathbb{D}^\varepsilon\,\mathbf{U}^{n}\,.\label{eq:matrix}
\end{equation}
This linear system is a rewriting of the relations \eqref{eq:h_meanfree}, \eqref{eq:scheme_lambda}, \eqref{eq:scheme_h_FP} in the Fokker-Planck case (or \eqref{eq:scheme_h_BGK} in the BGK case). The corresponding matrices $\mathbb{M}^\varepsilon_{\text{col}}$ are defined by
\[
 \mathbb{M}^\varepsilon_{\text{col}}\ =\ \left(\begin{matrix}
                              \mathbb{I}_{N} & \mathbb{M}^{\lambda, h}\\
                              \mathbb{M}^{h, \lambda} & \mathbb{M}^{h,h}_{\text{col}}\\
                              0 & \mathbb{N}^{h,h}                              
                            \end{matrix}\right)\,,\qquad 
                            \mathbb{D}^\varepsilon\ =\ \left(\begin{matrix}
                                                                \mathbb{I}_{N} & 0\\
								0 & \ve^2 \mathbb{I}_{2L}\\
								0&0
								\end{matrix}\right)\,.
\]

Here, we denoted by $\mathbb{I}_K$ the $K\times K$ identity matrix and the other matrices are defined as follows. The matrix $\mathbb{M}^{\lambda, h}$ has non-zero coefficients
\[
\mathbb{M}^{\lambda, h}_{\X_i, K_{i+1,j}}\ =\ \frac{v_j\,\M_j\,\Delta v_j\,\Delta t}{2\,\Delta x_i}\,,\qquad \mathbb{M}^{\lambda, h}_{\X_i, K_{i-1,j}}\ =\ -\frac{v_j\,\M_j\,\Delta v_j\,\Delta t}{2\,\Delta x_i}\,,\quad\forall (i,j)\in\I\times\J\,.
\]
The matrix $\mathbb{M}^{h, \lambda}$ has non-zero coefficients
\[
\mathbb{M}^{h, \lambda}_{\textcolor{black}{K_{i,j},\X_{i+1}}}\ =\ \frac{v_j\,\Delta t}{2\,\Delta x_i}\,,\qquad \mathbb{M}^{h, \lambda}_{\textcolor{black}{K_{i,j},\X_{i-1}}}\ =\ -\frac{v_j\,\Delta t}{2\,\Delta x_i}\,,\quad\forall (i,j)\in\I\times\J\,.
\]
The matrix $\mathbb{M}^{h, h}_{BGK}$ has non-zero coefficients
\begin{multline*}
\quad\forall (i,j)\in\I\times\J\,\quad \forall k\in\J\,,\  k\neq j\,,\qquad \mathbb{M}^{h, h}_{K_{ij}, K_{ij}}\ =\ \ve^2 + \Delta t\,,\\[1em]
\mathbb{M}^{h, h}_{K_{ij}, K_{i+1,j}}\ =\ \ve\,\frac{v_j\,\Delta t\,(1-\Delta v_j\,\M_j)}{2\,\Delta x_i}\,,\qquad
\mathbb{M}^{h, h}_{K_{ij}, K_{i-1,j}}\ =\ -\,\ve\,\frac{v_j\,\Delta t\,(1-\Delta v_j\,\M_j)}{2\,\Delta x_i}\,,\\[1em]
\mathbb{M}^{h, h}_{K_{ij}, K_{i+1,k}}\ =\ -\,\ve\,\frac{\textcolor{black}{v_k}\,\Delta t\,\Delta v_k\,\M_k}{2\,\Delta x_i}\,,\qquad
\mathbb{M}^{h, h}_{K_{ij}, K_{i-1,k}}\ =\ \ve\,\frac{\textcolor{black}{v_k}\,\Delta t\,\Delta v_k\,\M_k}{2\,\Delta x_i}\,.
\end{multline*}
The matrix $\mathbb{M}^{h, h}_{FP}$ has the same coefficients as $\mathbb{M}^{h, h}_{BGK}$  except for the following \textcolor{black}{(outside the boundary in velocity)}
\begin{multline*}
\quad\forall (i,j)\in\I\times\J\textcolor{black}{^*}\,, \,\qquad \mathbb{M}^{h, h}_{K_{ij}, K_{ij}}\ =\ \ve^2 + \frac{\Delta t\,}{\Delta v_j \textcolor{black}{\M_j}}\left(\frac{\M_{j+1/2}^*}{\Delta v_{j+1/2}}+\frac{\M_{j-1/2}^*}{\Delta v_{j-1/2}}\right)\,,\\[1em]
\mathbb{M}^{h, h}_{K_{ij}, K_{i,j+1}}\ =\ -\,\frac{\Delta t\,\M_{j+1/2}^*}{\Delta v_j\,\M_j\,\Delta v_{j+1/2}}\,,\qquad
\mathbb{M}^{h, h}_{K_{ij}, K_{i,j-1}}\ =\ -\,\frac{\Delta t\,\M_{j-1/2}^*}{\Delta v_j\,\M_j\,\Delta v_{j-1/2}}\,.
\end{multline*}
\textcolor{black}{At the boundary, $\forall i\in\I$,}
\begin{gather*}
\textcolor{black}{\mathbb{M}^{h, h}_{K_{i,-L+1}, K_{i,-L+1}}\ =\ \ve^2 + \frac{\Delta t\,}{\Delta v_{-L+1} \textcolor{black}{\M_{-L+1}}}\left(\frac{\M_{-L+3/2}^*}{\Delta v_{-L+3/2}}\right)\,,}\\[1em]
\textcolor{black}{\mathbb{M}^{h, h}_{K_{i,-L+1}, K_{i,-L+2}}\ =\ -\,\frac{\Delta t\,\M_{-L+3/2}^*}{\Delta v_{-L+1}\,\M_{-L+1}\,\Delta v_{-L+3/2}}\,,}\\[1em]
\textcolor{black}{\mathbb{M}^{h, h}_{K_{i,L}, K_{i,L}}\ =\ \ve^2 + \frac{\Delta t\,}{\Delta v_{L} \textcolor{black}{\M_{L}}}\left(\frac{\M_{L-1/2}^*}{\Delta v_{L-1/2}}\right)\,,} \qquad 
\textcolor{black}{\mathbb{M}^{h, h}_{K_{i,L}, K_{i,L-1}}\ =\ -\,\frac{\Delta t\,\M_{L-1/2}^*}{\Delta v_{L}\,\M_{L}\,\Delta v_{L-1/2}}\,.}
\end{gather*}
Finally the matrix $\mathbb{N}^{h, h}$, which comes from \eqref{eq:h_meanfree}, has non-zero coefficients
\[
  \mathbb{N}^{h, h}_{\X_{i}, K_{i,j}}\ =\ \M_j\,\Delta v_j,\qquad\forall (i,j)\in\I\times\J\,.
\]

\section{Numerical simulations}\label{sec:numres}

  In this final section, we shall present simulations for the overdetermined micro-macro versions $(\tilde{\mathbf{S}}_{\lambda,h}^\ve)$ of our numerical schemes. We will address numerically all the features of our numerical methods, namely their AP properties in a first subsection, and the discrete hypocoercivity in the next ones, along with some fine properties of the numerical solutions. 
  
  The numerical method has been written in\textit{ Python 3.6}, using \textit{Jupyter Lab}, and is available, together with all the test cases presented in this section in \cite{NumericalScheme:binder}.

    Unless stated otherwise, in all the simulations we will take $L = 20$, namely  $41$ points in a velocity box $[-v_*,v_*]$, with $v_* = 8$,  and $N = 51$ points in the torus $[0,1]$. Given that our schemes are unconditionally stable, we shall choose $\Delta t = 0.1$.
    
    \subsection{Numerical investigation of the AP Property}~

      \mysubsection{Test 1. Convergence toward the heat equation.} Let us start by investigating numerically the AP property of our schemes. We choose as an initial condition the following smooth, far-from-equilibrium distribution
      \begin{equation}
        \label{eq:initCondFar}
        f_0(x,v) := \frac{1}{\sqrt{2\pi}}\,v^4\,e^{-v^2/2}\,\frac{1+\cos(4\pi x)}{2}, \quad \forall x \in \mathbb{T}, v \in [-v_*,v_*].
      \end{equation}
      
\begin{figure}
        \begin{center}
        \input{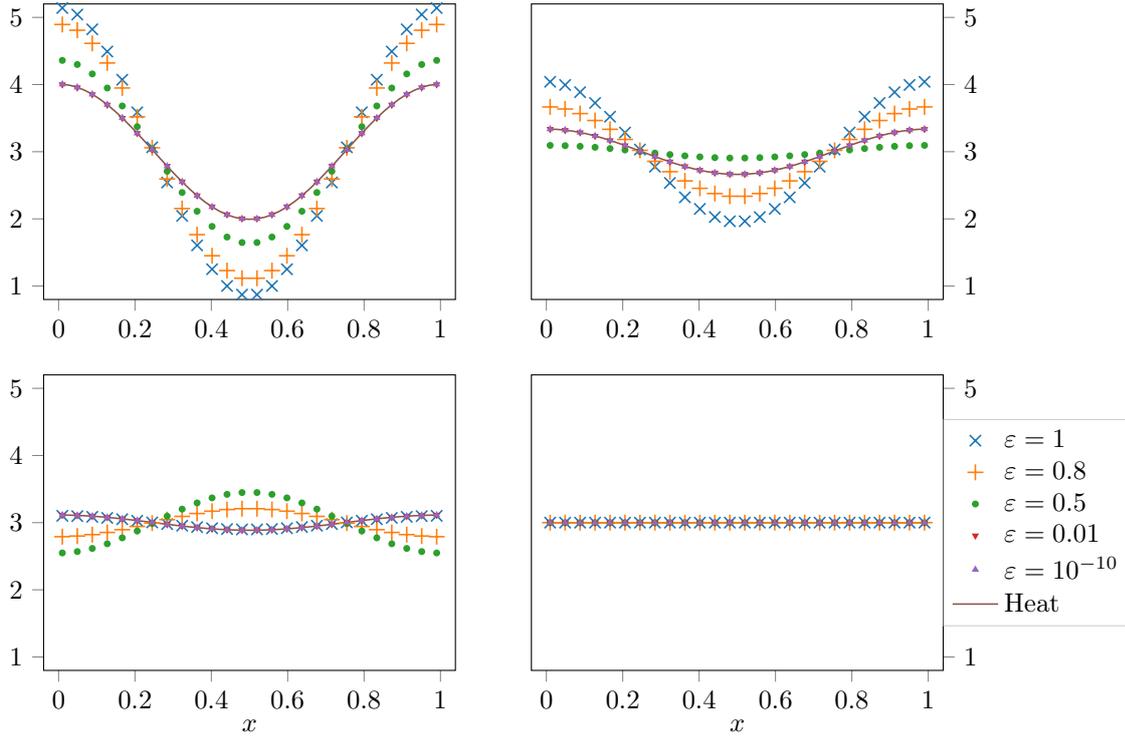}    
        \end{center}  
        \caption{\textbf{Test 1. }Comparison of the solution to the heat equation obtained from the scheme \eqref{eq:heat_discrete} (solid line) with the solution obtained with the \textbf{Fokker-Planck} scheme \eqref{eq:scheme_h_FP} with different $\ve$, at times $t=0.05$, $0.1$, $0.15$ and $10$.}
        \label{fig:HeatvsKineticFP}
      \end{figure}

      We solve first the kinetic Fokker-Planck equation with scheme \eqref{eq:scheme_h_FP}, for \Mar{different values of $\ve$}, and compare the discrete density $\rho^\ve$ at different times with the solution to the discrete heat equation given by \eqref{eq:heat_discrete}. The results can be found in Figure \ref{fig:HeatvsKineticFP}, where we observe a good agreement between the convergence result of Theorem \ref{th:AP} and our simulations.

\begin{figure}
        \begin{center}
        \input{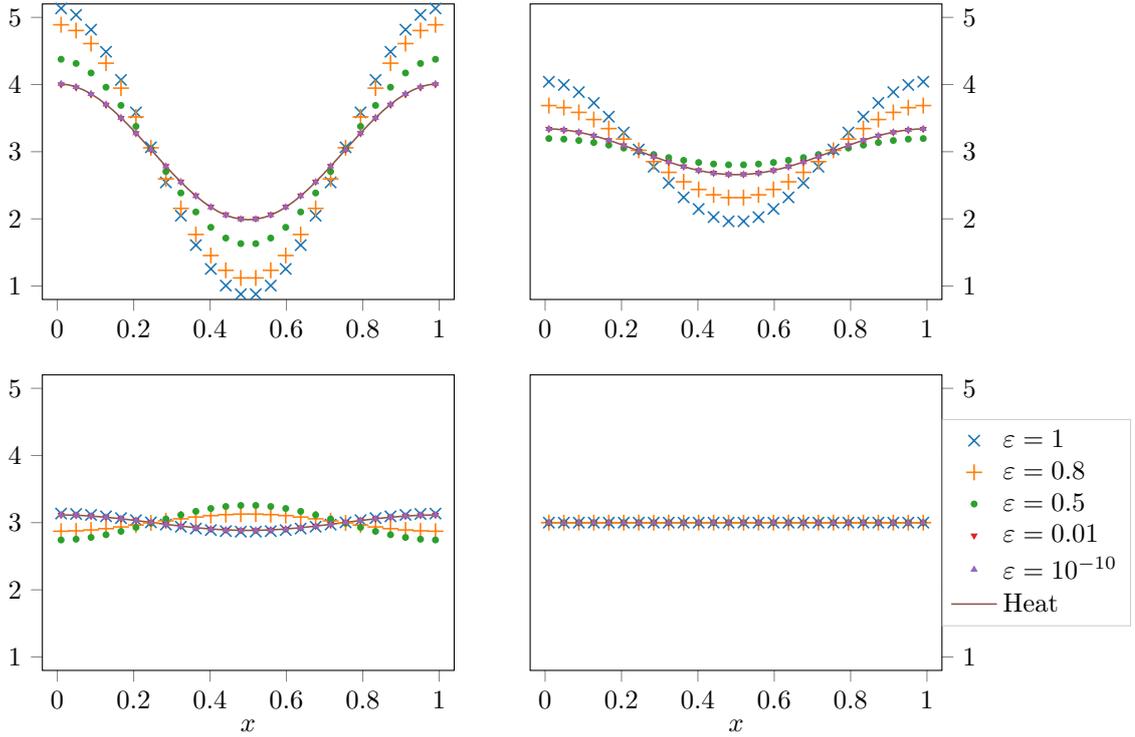}    
        \end{center}  
        \caption{\textbf{Test 1. }Comparison of the solution to the heat equation obtained from the scheme \eqref{eq:heat_discrete} (solid line) with the solution obtained with the \textbf{BGK} scheme \eqref{eq:scheme_h_BGK} with different $\ve$, at times $t=0.05$, $0.1$, $0.15$ and $10$.}
        \label{fig:HeatvsKineticBGK}
      \end{figure}

       We then perform the same numerical investigations for the linear BGK equation with scheme \eqref{eq:scheme_h_BGK}. Again, we observe in Figure \ref{fig:HeatvsKineticBGK} the convergence with respect to $\ve$ of the density $\rho^\ve$ towards the solution to the discrete heat equation $\rho$, as predicted by Theorem \ref{th:AP}. 

      We will now investigate the hypocoercivity properties of our schemes, with various numerical experiments in the Fokker-Planck and BGK cases.

    \subsection{The Fokker-Planck case}~

      \mysubsection{Test 2. Trend to equilibrium.}
        We shall fix $\ve = 1$ for this numerical test case. 
        
        Let us first start to verify the results of Theorem \ref{th:hypocoercivity_discrete}, namely that the discrete solution $f^\ve$ to the kinetic Fokker-Planck equation converges towards $\mu_f \M$ in the discrete, weighted, $L^2$ topology, at an exponential rate.
        
\begin{figure}
          \begin{center}
            \begin{tabular}{cc}
              \includegraphics[scale=1]{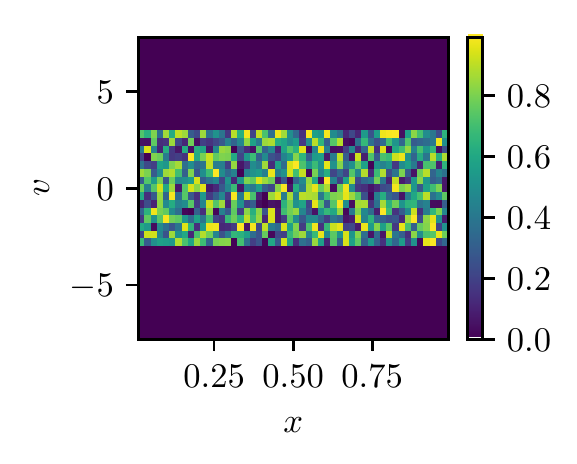} & \input{./Pics/TrendToEq_Equationkinetic_FP_InitDat2_Eps1} 
            \end{tabular}
          \end{center}
          \caption{\textbf{Test 2.} Left: Initial data in the $(x,v)-$phase-plane. Right: Time evolution of the weighted $L^2$ norm of the difference between $f^\ve$ and the global equilibrium, $\ve =1$ in the \textbf{Fokker-Planck} case, and observed exponential rate.}
          \label{fig:Trend2Eq_RandomTruncat}
        \end{figure}

        We first perform this experiment with an initial data chosen to be random, uniformly distributed on the $(x,v)-$plane, but supported in $v \in [-3,3]$ (also shown in the left side of Figure~\ref{fig:Trend2Eq_RandomTruncat}).
        We observe in the rightmost part of Figure \ref{fig:Trend2Eq_RandomTruncat} that the numerical rate of convergence of the discrete norm $\|f - \mu_f \M \|_{2,\gamma}$ obtained with the scheme \eqref{eq:scheme_h_FP} is indeed exponential. It also matches the rate of convergence in the spatial $L^2$ norm of the density $\rho^\ve$ towards the global mass $\mu_f$. The numerical hypocoercivity rate for this test case is $1.003$, which is close to the rate observed\footnote{\Tom{Note that the quantity we compute in our paper is the square root of the one computed in \cite{Dujardin}, explaining the factor 2 between the observations.}} in the paper~\cite{Dujardin}.
        
\begin{figure}
          \begin{center}
            \begin{tabular}{cc}
              \includegraphics[scale=1]{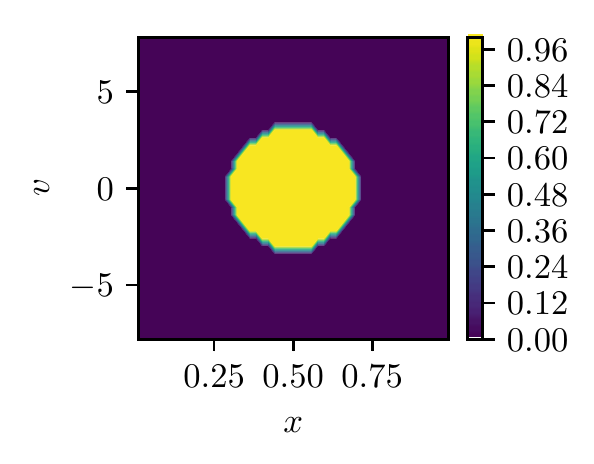} & \input{./Pics/TrendToEq_Equationkinetic_FP_InitDat3_Eps1} 
            \end{tabular}
          \end{center}
          \caption{\textbf{Test 2.} Left: Initial data in the $(x,v)-$phase-plane. Right: Time evolution of the weighted $L^2$ norm of the difference between $f^\ve$ and the global equilibrium, $\ve =1$, \textbf{Fokker-Planck} case.}
          \label{fig:Trend2Eq_Ball}
        \end{figure}

These observations are confirmed by choosing another initial condition, such as  the indicator of a ball in the $(x,v)-$phase plane, as seen on the left of Figure \ref{fig:Trend2Eq_Ball}.
        We observe in the \Tom{rightmost} part of that figure the same behavior than in the previous case. The numerical rate of decay is \Tom{now  $1.96$,} quite similar to the value observed in \cite{Dujardin}.
        
       %Note that the exact rate in the space homogeneous case is exactly $1$ (see again \cite{Dujardin} for a short proof). The optimal rate in the full space inhomogeneous case is not known, \Tom{but seems to be strongly dependent on the initial condition}.
        
\begin{figure}
        \begin{center}
        \includegraphics[scale=1]{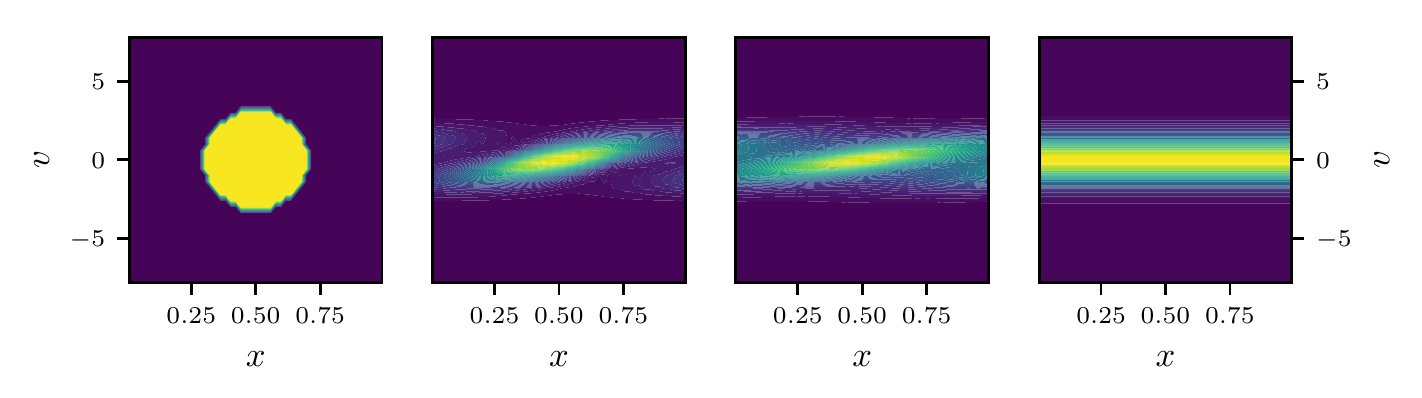}
        \end{center}
        \caption{\textbf{Test 2.} Snapshot of the particle distribution function $f^\ve(t,x,v)$ in the $(x,v)-$phase-plane, at times $t = 0$, $0.3$, $0.6$, and $30$.}
        \label{fig:SnapshotFP}
      \end{figure}

        Finally, we also present in Figure \ref{fig:SnapshotFP} a snapshot of the particle distribution function  in the $(x,v)-$phase-plane at different times, illustrating its convergence towards the global Maxwellian distribution.
        
      \mysubsection{Test 3. Dependency of the exponential rate on $\ve$.} For this final Fokker-Planck test, we  study the influence of $\ve$ on the rate  of exponential decay of the quantity $\|f^\ve - \mu_f \M \|_{2,\gamma}$.         
        We choose initially a uniformly distributed distribution $f_0$ in the $(x,v)-$plane, with values between 0 and 1, and we let $\ve$ vary from the purely kinetic regime $\ve = 1$ to the purely diffusive regime $\ve =0$, without ever changing the other grid parameters. 
        
\begin{figure}
          \begin{center}
            \includegraphics{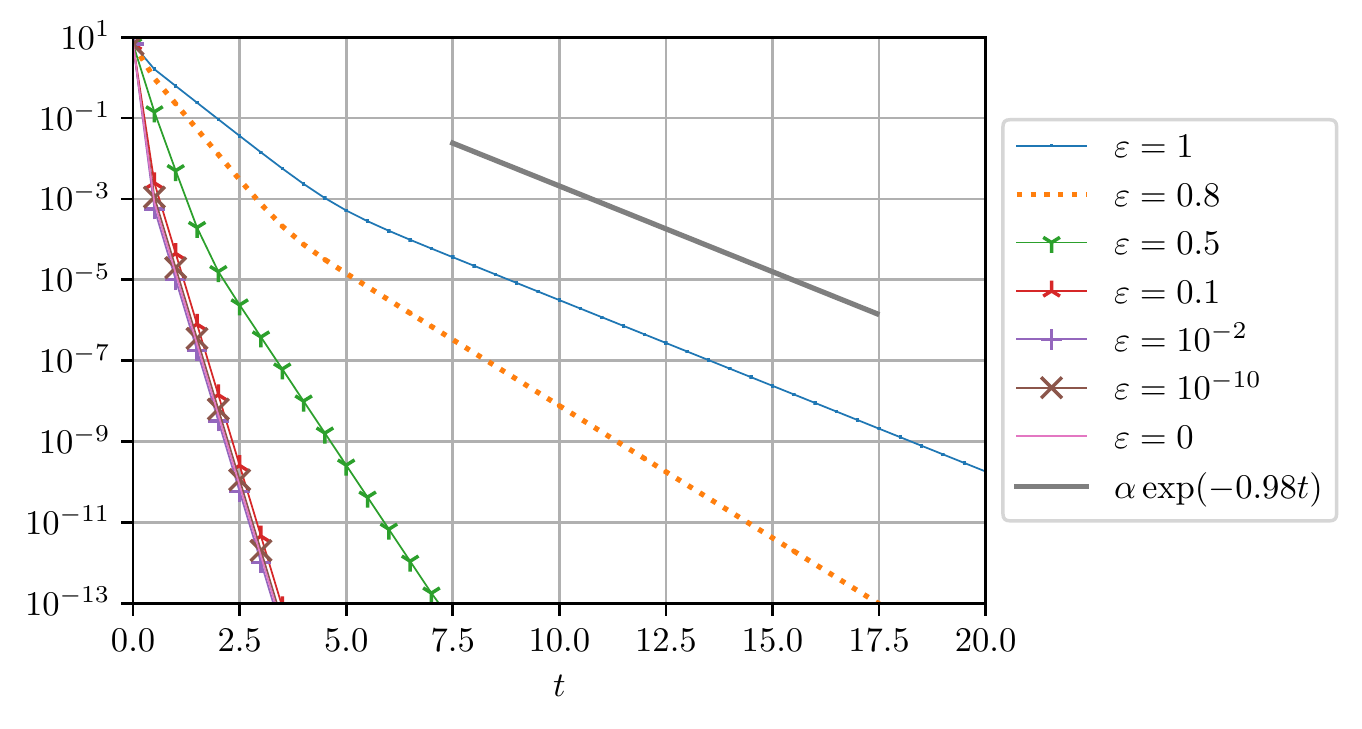}
          \end{center}
          \caption{\textbf{Test 3.} Comparison of the rate of convergence of $\|f^\ve - \mu_f \M \|_{2,\gamma}$ for different values of $\ve$, in the \textbf{Fokker-Planck} case.}
          \label{fig:CompEpsFP}
        \end{figure}

        We represent the results of this numerical investigation in Figure \ref{fig:CompEpsFP}, by plotting the different numerical values of $\|f^\ve(t) - \mu_f \M \|_{2,\gamma}$.
        We observe that the rate of exponential decay seems to have a lower bound which is uniform  on $\ve$, illustrating the fact that the result of Theorem \ref{th:hypocoercivity_discrete} is suboptimal concerning the rate $\beta$ obtained. This could of course be improved by tracking exactly the dependence in $\ve$ in the estimates used in the proof of Theorem~\ref{th:hypocoercivity_discrete}.
        
        To make this latter statement more explicit, we show in Table \ref{tab:rates} the value of the rate of decay, as a function of $\ve$.  This rate decreases with $\ve$ towards a limit value.
%        1.e+00, 8.e-01, 5.e-01, 1.e-01, 1.e-02, 1.e-10, 0.e+00
%        -0.9761760909710392,
%  -1.5044695333438702,
%  -3.647457811290491,
%  -8.044263040178471,
%  -8.047078932532276,
%  -8.04710740805842,
%  -8.047107120360497
      \begin{table}
        \begin{center}
        \begin{tabular}{|c|c|}
          \hline
          $\ve$ & Rate \\ \hline
          $1$         & -0.98 \\ \hline
          $0.8$      & -1.5 \\ \hline
          $0.5$      & -3.65 \\ \hline
          $10^{-1}$  & -8.04 \\ \hline
          $10^{-2}$  & -8.05 \\ \hline
          $10^{-10}$ & -8.05 \\ \hline
          0          & -8.05 \\ \hline
        \end{tabular}   
        \end{center}
        \caption{\textbf{Test 3.} Rate of exponential decay of $\|f^\ve(t) - \mu_f \M \|_{2,\gamma}$ with respect to $\ve$.}
        \label{tab:rates}       
      \end{table}      
      
    \subsection{The BGK case}~
    
      \mysubsection{Test 4. Oscillations during the macroscopic relaxation to equilibrium.}
        We now \Mar{consider}  the question of the macroscopic relaxation towards equilibrium. So far, we investigated the behavior of the rate of relaxation in $L^2(\dD x\dD\gamma)$ of $f^\ve$ towards its global equilibrium $\mu_f \M$, namely kinetic relaxation. We will now turn to a numerical study of the rate of relaxation in $L^2(\dD x)$ of the velocity averages of $f^\ve$, namely $\rho^\ve$, towards their constant-in-space global equilibrium $\mu_f$. We shall call this later case \emph{macroscopic} relaxation.
        
        The behavior of this macroscopic quantity has already been presented in Figures \ref{fig:Trend2Eq_RandomTruncat} and \ref{fig:Trend2Eq_Ball}, where one can notice that the macroscopic relaxation seems bounded by the kinetic one. Nevertheless, this macroscopic behavior is actually richer, and can even be not monotonous.
        To study it more precisely, we shall consider the following close-to-equilibrium initial condition
        \begin{equation}
          \label{eq:initCondClose}
           f_0(x,v) :=\left (1+\cos\left (\frac{2\pi x}{R}\right )\right ) \M(v), \quad \forall x \in [0,R], \quad v \in [-v_*,v_*].
        \end{equation}
        Our goal is to change the size $R$ of the torus, to see how it affects the behavior of the macroscopic relaxation. Note that in the previous test cases, $R$ was equal to $1$. 
        
\begin{figure}
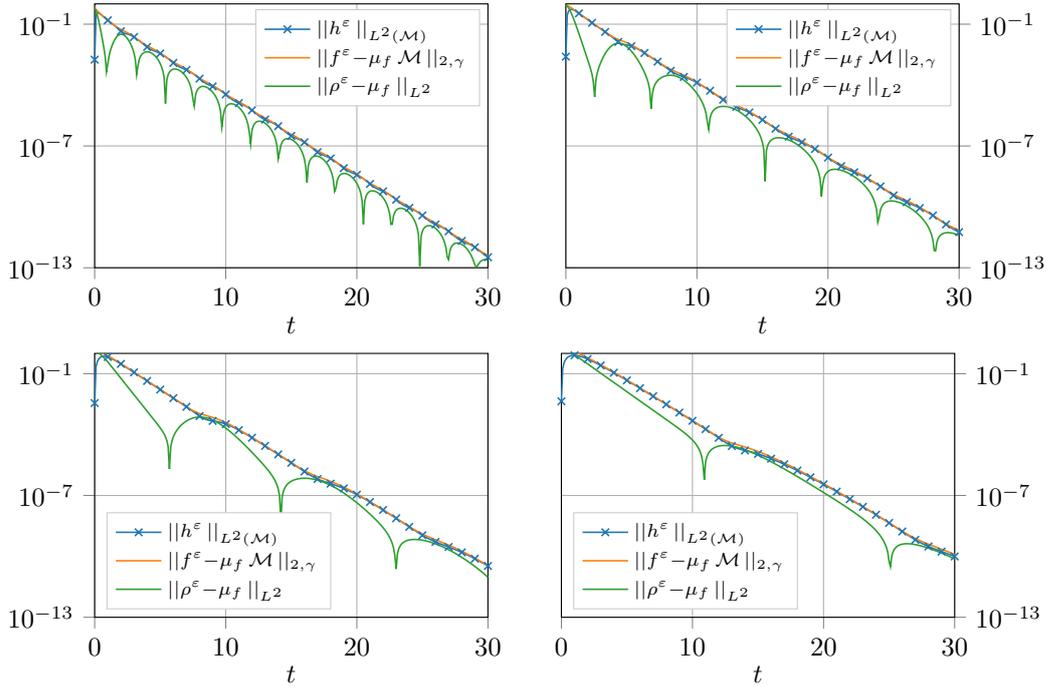

          \begin{center}
          \begin{tabular}{cc}          
            \input{Pics/TrendToEq_EquationBGK_InitDat4_Eps1_X025pi} &
            \input{Pics/TrendToEq_EquationBGK_InitDat4_Eps1_X05pi}\\
            \input{Pics/TrendToEq_EquationBGK_InitDat4_Eps1_Xpi} & \input{Pics/TrendToEq_EquationBGK_InitDat4_Eps1_X15pi}
          \end{tabular}
          \end{center}
          \caption{\textbf{Test 4.} Influence of the torus size $R$ on the macroscopic relaxation, for $R = \pi/4$, $\pi/2$, $\pi$, and $3\pi/2$, in the \textbf{BGK} case.}
          \label{fig:Oscillations}
        \end{figure}

        Figure \ref{fig:Oscillations} presents, for the BGK case, the time evolution of the local relaxation $\|f^\ve(t) - \mu_f \M \|_{2,\gamma} $, and of the macroscopic relaxation $ \| \rho^\ve(t) - \mu_f \|_{L^2(\dD x)}$, for $4$ box sizes $R$. The $L^2(\M\dD x \dD v)$ norm of the scaled, microscopic unknown $h^\ve:= (f^\ve - \mu_f \M)/\M$   is also shown on these figures.
        The scaling constant $\ve$ is chosen equal to $1$ in all these simulations.
        
        We observe that, as expected through Theorem \ref{th:hypocoercivity_discrete}, the local relaxation rate is monotone and exponential.  This rate is inversely proportional to the box size $R$.
        Nevertheless, the macroscopic rate of relaxation is not monotone, exhibiting exponentially damped  oscillations.
        We also observe that the upper envelope of these oscillations seems to be given by 
        $\|h^\ve(t)\|_{L^2(\M\dD x \dD v)}$, and that this quantity seems to decay in a monotone but slightly oscillating way.

       Inspired by the analysis conducted in \cite{FiMoPa:2006}, we will give a simple interpretation for \Mar{these} oscillatory behaviors. 
       Since the BGK equation we consider is linear, the initial value \eqref{eq:initCondClose} we chose corresponds to observing the time evolution of the Fourier mode in $x$ of frequency $k = 2\pi/R$.
       Going to Fourier in $x$ in the BGK equation \eqref{eq:Collision} amounts to solve the following integral equation
        \begin{equation}
          \label{eq:gFourier}
          \partial_t g = \left(\int_\R g(v) \,     \dD v\right) \, \M(v) - \left (1+ i k v\right ) g.
        \end{equation}
        
        The right hand side of \eqref{eq:gFourier} is the sum of an operator of multiplication by $\mu(v) := - \left (1+i k v\right )$ and of a rank-1 projector (towards the global Maxwellian), that we shall denote by $\mathcal K$. 
        This operator then generates a strongly continuous semigroup $T_k(t)$. 
        Moreover, according to Weyl's Theorem \cite{kato:1966}, its spectrum in $L^2$ contains the spectrum of the multiplication operator, namely the numerical range of $\mu$: 
        \[
            \{-1 + i y: y \in \R\}\, ,
        \]       
        together with some discrete eigenvalues located in the complex plane on the right of this set.
        Since we proved in Theorem \ref{th:hypocoercivity_discrete} that the linear BGK equation is hypocoercive with a rate at least $\beta/2$, necessarily, such a discrete nonzero eigenvalue $\lambda \in \CC^*$ verify
        \[
            -1 \leq \Re e \left ( \lambda \right ) \leq -\beta/2.
        \]
        Let $g_\lambda \in L^2(\dD x \dD \gamma)\setminus \{0\}$ be the associated eigenvector. One then has
        \begin{equation}
            \label{eq:eigProb}
            \mathcal K g_\lambda = \left (1 + \lambda + ikv \right)g_\lambda.
        \end{equation}
        Taking the $L^2(\dD x \dD \gamma)-$inner product of this relation with $g_\lambda$ yields
        \begin{align}
            \lambda  & = -1 + \frac{\langle \mathcal K g_\lambda, g_\lambda \rangle_{L^2(\dD x \dD \gamma)}}{\| g_\lambda \|_{L^2(\dD x \dD \gamma)}^2} - ik \frac{\int_\RR v \, |g_\lambda|^2 \,\dD \gamma}{\| g_\lambda \|_{L^2(\dD x \dD \gamma)}^2} \notag \\
            & = -\left[ 1 - \frac{\left( \int_\RR g_\lambda (v) \, \dD v \right )^2}{\| g_\lambda \|_{L^2(\dD x \dD \gamma)}^2} \right] - ik \frac{\int_\RR v \, |g_\lambda|^2 \,\dD \gamma}{\| g_\lambda \|_{L^2(\dD x \dD \gamma)}^2}. \label{eq:lambdaglambda}
        \end{align}
        
        Note that the case of the space-homogeneous solution to \eqref{eq:Collision} corresponds to $k = 0$, so that taking the integral of \eqref{eq:eigProb} yields that $g_\lambda$ is of zero mass. Hence, in that particular case, \eqref{eq:lambdaglambda} implies that $\lambda =-1$ is the only nonzero eigenvalue.%, a result which was already proved in \cite{Dujardin}.
        
        \Mar{Finally}, we have given a formal proof that all the nonzero eigenvalues of  \eqref{eq:gFourier} can be written as $\lambda = -\kappa + i k u$, with $0<\kappa<1$ and $u \in \RR$. 
        The $L^2(\dD x \dD \gamma)$ norm of the solutions to this equation will then behave as
        \[ \Re e\exp( (-\kappa+ik u)t) = \exp(-\kappa t)\cos(k\,u\,t), \]
        namely decay exponentially in time, while oscillating at a period proportional to $k = 2\pi/R$.
        This is also in good agreement with the results from \cite{lafitte2016high}, where the spectrum of a fully discretized linear BGK equation is computed and shown to behave similarly.
        
        To verify this interpretation experimentally, we compute in Table \ref{tab:Oscillations} the oscillation frequency $\nu$ with respect to the box size $R$. We observe that the product $R\nu$ appears indeed almost constant, yielding $u \simeq 0.35$ in the ongoing interpretation and giving more evidences that our numerical scheme preserves expected theoretical behaviors of the continuous equation.
        
        \begin{table}
          \begin{center}
          \begin{tabular}{|c|c||c|}
            \hline
            Torus length $R$ &  Oscillation period $\nu^{-1}$ & Speed $R\nu$ \\ \hline
            $\pi /4$         &            2.31                & 0.340       \\ \hline
            $\pi /2$         &            4.33                 & 0.363       \\ \hline
            $\pi   $         &            8.67                 & 0.362       \\ \hline
            $3\pi/2$         &            13.5                 & 0.350       \\ \hline
          \end{tabular}
          \end{center}
          \caption{\textbf{Test 4.} Influence of the torus length on the oscillation period of $\|\rho^\ve - \mu_f\|_{L^2}$, in the \textbf{BGK} case.}
          \label{tab:Oscillations}
        \end{table}

    \subsection{\Tom{The Non-Gaussian BGK case}}\label{s:nongauss}~
    
        \mysubsection{\Tom{Test 5. Relaxation towards a non-Gaussian equilibrium.}} \Tom{As a final test case, we consider the relaxation of the solution to the linear BGK equation towards the following equilibrium distribution with polynomial \Mar{tail}:
        \begin{equation}
            \label{eq:nongaussianEQ}
            M(v) \, = \, \frac{1}{1+0.1|v|^6} \left(\cos(\pi v)+1.1\right).
        \end{equation}
        This distribution verifies the hypothesis \eqref{hyp:MaxwBGK}, namely that once projected on the velocity grid, it is of discrete mass 1, and has uniformly bounded discrete moments of order 2 and 4. According to Theorem \ref{th:hypocoercivity_discrete}, the solution to the discrete BGK equation solved with the scheme \eqref{eq:scheme_f_BGK} will exhibit a discrete hypocoercive behavior.}
        
\begin{figure}
          \begin{center}
            \begin{tabular}{cc}
              \includegraphics[scale=1]{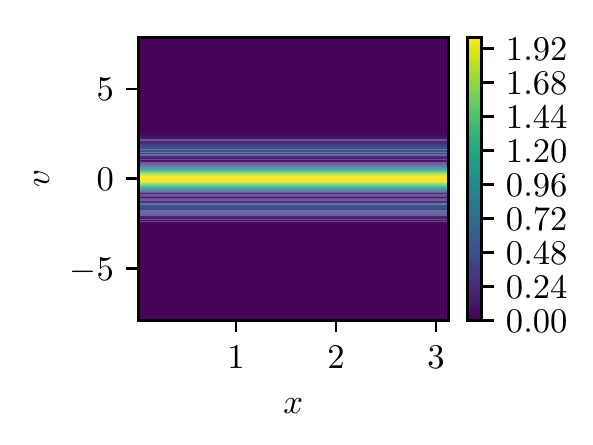} & \input{./Pics/TrendToEq_EquationNonGaussianBGK_InitDat5_Eps1} 
            \end{tabular}
          \end{center}
          \caption{\Tom{\textbf{Test 5.} Left: Equilibrium distribution ($t=40$) in the $(x,v)-$phase-plane. Right: Time evolution of the weighted $L^2$ norm of the difference between $f^\ve$ and the global equilibrium, $\ve =1$, \textbf{Non-Gaussian BGK} case.}}
          \label{fig:Trend2Eq_nongaussian}
        \end{figure}

        \Tom{We choose the very far from equilibrium \Mar{initial condition} \eqref{eq:initCondFar} \Mar{represented} on the left side of Figure \ref{fig:Trend2Eq_nongaussian}. Because of the fine structures in the equilibrium profile, we chose a more refined grid than in the previous test cases, namely $L=35$, $N=101$ and $\Delta t = 0.01$.
        We observe on the right of Figure \ref{fig:Trend2Eq_nongaussian} that the solution \Max{decays exponentially} towards the non-Gaussian equilibrium \eqref{eq:nongaussianEQ}.}

\section{Appendix. On discrete Poincaré inequalities}

For the sake of completeness, we give in this Appendix the proofs of Poincaré inequalities we used in this article. We first start with the discrete Gaussian Poincaré inequality on bounded domain.

\begin{proof}[Proof of Lemma \ref{lem_poincare_v}]
 We follow the same guidelines as in \cite[Proposition 4.7]{Dujardin}. Replacing if necessary $f$ by $f-\rho$, where $\rho=\sum_{j\in\J}f_j\Delta v_j$, we can assume without loss of generality that $\rho=0$. First we can remark that
 \begin{align*}
     \|f\|_{2,\gamma}^2&=\sum_{j\in\J}(f_j\gamma_j)^2\M_j\Delta v_j=\frac{1}{2}\sum_{(j,k)\in\J^2}(f_k\gamma_k-f_j\gamma_j)^2\M_j\M_k\Delta v_j \Delta v_k\\
     &= \sum_{j<k}(f_k\gamma_k-f_j\gamma_j)^2\M_j\M_k\Delta v_j\Delta v_k.
 \end{align*}
 For $j<k$, we have
 \[
 f_k\gamma_k-f_j\gamma_j=\sum_{l=j}^{k-1}(f_{l+1}\gamma_{l+1}-f_l\gamma_l).
 \]
 Then, we obtain using Cauchy-Schwarz inequality 
 \begin{align*}
     \|f\|_{2,\gamma}^2&=\sum_{j<k}\left(\sum_{l=j}^{k-1}\frac{f_{l+1}\gamma_{l+1}-f_l\gamma_l}{\sqrt{\Delta v_{l+\frac{1}{2}}}}\sqrt{\Delta v_{l+\frac{1}{2}}} \right)^2\M_j\M_k\Delta v_j\Delta v_k\\
     &\leq \sum_{j<k}\left(\sum_{l=j}^{k-1}\left( D_v(f\gamma)_{l+\frac{1}{2}}\right)^2\Delta v_{l+\frac{1}{2}} \right)(v_k-v_j)\M_j\M_k\Delta v_j\Delta v_k.
 \end{align*}
 Let us introduce $F=(F_j)_{j\in\J}$ defined by
 \begin{align*}
     &F_j:=-\sum_{l=j}^{-1}\left( D_v(f\gamma)_{l+\frac{1}{2}}\right)^2\Delta v_{l+\frac{1}{2}}\quad \text{ if }j\leq -1,\\
     &F_j:=\sum_{l=0}^{j-1}\left( D_v(f\gamma)_{l+\frac{1}{2}}\right)^2\Delta v_{l+\frac{1}{2}} \quad\text{ if }j\geq 0,
 \end{align*}
 so that for $j<k$, we have 
 \[
 F_k-F_j=\sum_{l=j}^{k-1}\left( D_v(f\gamma)_{l+\frac{1}{2}}\right)^2\Delta v_{l+\frac{1}{2}}.
 \]
 It yields
 \begin{align*}
     \|f\|_{2,\gamma}^2&\leq \sum_{j<k}(F_k-F_j)(v_k-v_j)\M_j\M_k\Delta v_j\Delta v_k\\
     &\leq \frac{1}{2}\sum_{(j,k)\in\J^2}(F_k-F_j)(v_k-v_j)\M_j\M_k\Delta v_j\Delta v_k
 \end{align*}
 since $(F_k-F_j)(v_k-v_j)=(F_j-F_k)(v_j-v_k)$. 
 
 Developing the right-hand side and using \eqref{hyp:MaxwFP}, we get
 \[
 \|f\|_{2,\gamma}^2\leq \sum_{j\in\J}F_jv_j\M_j\Delta v_j.
 \]
 Now, using the definition of $\M_j$ given in \eqref{hyp:MaxwFP}, we obtain
 \[
 \|f\|_{2,\gamma}^2\leq \sum_{j\in\J}F_j(\M^*_{j-\frac{1}{2}}-\M^*_{j+\frac{1}{2}}).
 \]
 Applying a discrete integration by parts, it gives
 \[
 \|f\|_{2,\gamma}^2\leq\sum_{j\in\Jhalf}(F_{j+1}-F_j)\M^*_{j+\frac{1}{2}},
 \]
 which finally yields the result since by definition,
 \[
 F_{j+1}-F_j=\left( D_v(f\gamma)_{j+\frac{1}{2}} \right)^2\Delta v_{j+\frac{1}{2}}=\left( D_v\left(\frac{f}{\M}\right)_{j+\frac{1}{2}} \right)^2\Delta v_{j+\frac{1}{2}}.
 \]

\end{proof}

\color{black}
Our numerical scheme and every result in the paper generalize in higher dimension on Cartesian meshes by tensorization. Hereafter we just show how it is done for one important point of our analysis: the discrete Gaussian Poincaré inequality.
\begin{cor}[Discrete Gaussian Poincaré inequality in higher dimension]\label{lem_poincare_v_multid}
Let $d\geq1$. Given a discrete unidimensional Maxwellian $(\M_{k})_{k\in\J}$ and $(\M_{k+1/2}^*)_{k\in\J^*}$ satisfying \eqref{hyp:MaxwFP}, a multi-index $j = (j_1,\dots,j_d)$ and $e\in\{1,\dots,d\}$, let us define multidimensional Maxwellians at cell centers and edge/face centers respectively by
\[
\M_j\ :=\ \prod_{e=1}^d\M_{j_e}\,,\quad \M_{j+\delta_e/2}\ :=\ \M_{j_{e}+1/2}\prod_{\substack{l=1\\l\neq e}}^d\M_{j_{l}}\,,
\]
\Mar{where $\delta_e$ denotes the multi-index whose all components are zero, except the $e$-th.}
We also write
$
\Delta v_j\ :=\ \prod_{e=1}^d\Delta v_{j_e}$ and $\Delta v_{j+\delta_e/2}\ :=\ \Delta v_{j_{e}+1/2}\prod_{\substack{l=1\\l\neq e}}^d\Delta v_{j_{l}}
$. Then for all $f = (f_{j})_{j\in\J^d}\in\RR^{\J^d}$,  one has that
\begin{multline}\label{ineq_Poincare_v_multid}
\sum_{j\in\J^d}|f_j - \rho\M_j|^2\,\M_j^{-1}\,\Delta v_j\ \leq\ \\\sum_{e=1}^d\sum_{j\in\J^{e-1}\times\J^*\times\J^{d-e}}\left(\frac{\Mar{(f/\M)_{(j_1,\dots,j_e+1,\dots,j_d)} - (f/\M)_{(j_1,\dots,j_e,\dots,j_d)}}}{\Delta v_{j+\delta_e/2}}\right)^2\Delta v_{j+\delta_e/2}\,\M_{j+\delta_e/2}\,,
\end{multline}
where $\rho = \sum_{j\in\J^d} f_j \Delta v_j$.
\end{cor}
\begin{proof}
\Mar{Let us introduce the notation $f^{(e)}_j = \sum_{(k_1,\dots,k_e)\in\J^e} g_{(k_1,\dots,k_e,j_{e+1},\dots,j_d)}\prod_{l=1}^e\Delta v_{k_l}$.} Observe that $f^{(0)} = f$ and that   $f^{(d)} = \rho$. Then, use Pythagoras equality and \Mar{apply $e-1$ times Jensen inequality }(since $\sum_k\M_k\Delta v_k = 1$) to get
\[
\begin{array}{rl}
&\sum_{j\in\J^d}|f_j - \rho\M_j|^2\,\M_j^{-1}\,\Delta v_j\\[.5em]
=& \sum_{j\in\J^d}|\sum_{e=1}^d(f_j^{(e-1)} \prod_{u=1}^{e-1}\M_{j_u} - f_j^{(e)} \prod_{u=1}^{e}\M_{j_u})|^2\,\M_j^{-1}\,\Delta v_j\\[.5em]
=& \sum_{e=1}^d\sum_{j\in\J^d}|f_j^{(e-1)} \prod_{u=1}^{e-1}\M_{j_u} - f_j^{(e)} \prod_{u=1}^{e}\M_{j_u}|^2\,\M_j^{-1}\,\Delta v_j\\[.5em]
=& \sum_{e=1}^d\sum_{(j_e,\dots,j_d)\in\Mar{\J^{(d-e+1)}}}|f_j^{(e-1)} - f_j^{(e)}\M_{j_e}|^2\,\prod_{u=e}^{d}\M_{j_u}^{-1}\Delta v_{j_u}\\[.5em]
\leq&\sum_{e=1}^d\sum_{j\in\J^{d}}|f_j - \M_{j_e}\sum_{k_e\in\J}f_{(j_1,\dots,j_{e-1},k_e,j_{e+1},\dots,j_d)}\Delta v_{k_e}|^2\,\M_{j}^{-1}\,\Delta v_{j}
\end{array}
\]
and one concludes by using $d$ times the 1D discrete Gaussian Poincaré inequality \eqref{ineq_Poincare_v}.
\end{proof}
\color{black}
We \Max{finally} provide a proof for the discrete Poincaré inequality on the (spatial) torus, \Max{in 1D and on a uniform mesh to simplify the presentation. The tensorization of the inequality may be done as before if one wants to generalize it to the multidimensional Cartesian setting.}

\begin{proof}[Proof of Lemma \ref{lem:poincare_discrete}]
  Let us assume that the grid is uniform in $x$, namely that $\Delta x_i = \Delta x >0$, for any $i \in \mathcal I = \{0,\ldots, N-1\}$. We will also assume that the indices are $N-$periodic. 
  This will allow to write a much simpler proof, using discrete Fourier transform (DFT). The \Mar{result} yet holds true in the general case, as long as $N$ is \textbf{odd}.
  
  Let us define the DFT of an $N-$periodic vector $g \in \mathbb C^N$, denoted by $\mathcal{F}_N g \in \mathbb C^N$, by:
  \[
    \left (\mathcal{F}_N g\right )_k := \Delta x\sum_{l=0}^{N-1} g_l \, e^{-2 i \pi k l \Delta x}, \quad \forall k \in \{0,\ldots,N-1\}.
  \] 
  Such a discrete transform \Mar{verifies} the following discrete Parseval identity
  \begin{equation}
    \label{eq:DiscreteParseval}
    \sum_{l =0}^{N-1} | g_l|^2 \Delta x = \sum_{k = 0}^{N-1} \left | \left ( \mathcal F_N g \right )_k \right |^2 .
  \end{equation}
  Since $\phi$ is of zero mass, one has that $\left (\mathcal{F}_N \phi \right )_0 = 0$. Then using identity \eqref{eq:DiscreteParseval}, we have
  \begin{align}
      \| D_x \phi \|_{L^2}^2 & = \sum_{l =0}^{N-1} \left | (D_x \phi)_l \right |^2 \Delta x \notag \\
      & = \sum_{k = 0}^{N-1} \left | \left ( \mathcal F_N D_x \,\phi \right )_k \right|^2\,. \label{eq:L2D_xFn}
  \end{align}
  Since $(D_x  \phi)_l = (\phi_{l+1} - \phi_{l-1})/(2 \Delta x)$, one has using the definition of the DFT that for any $k \in \{1,\ldots,N-1\}$
  \begin{align}
    \left ( \mathcal F_N D_x \, \phi \right )_k & = -\frac{e^{2 i \pi k \Delta x} - e^{-2 i \pi k \Delta x}}{2 \Delta x} \left ( \mathcal F_N \, \phi\right )_k \notag \\
    & = -i\, \frac{\sin (2 \pi k \Delta x)}{\Delta x} \left ( \mathcal F_N \phi\right)_k.    \label{eq:FnDxPhiShift} 
  \end{align} 
  If $N$ is even, then for $k=N/2$ the right hand side of \eqref{eq:FnDxPhiShift} will be equal to $0$. Since $N$ is odd, and $\left (\mathcal{F}_N \phi \right )_0 = 0$, one has
  \begin{equation}
     \| D_x \phi \|_{L^2}^2  \geq S_*^2 \, \sum_{k = 1}^{N -1} \left | \left ( \mathcal F_N \phi \right)_k \right |^2 
     \label{eq:almostPoincare}
  \end{equation}
  where 
  \begin{align}
      S_* & := \min_{k \in \{1,\ldots,N-1 \} } \left \{ \left | \frac{\sin (2 \pi k \Delta x)}{\Delta x} \right |  \right \}\notag \\
       & = \min \left \{\left | \frac{\sin (2 \pi \Delta x)}{\Delta x} \right | , \left |\frac{\sin (\pi (N-1) \Delta x)}{\Delta x} \right | \right \}
       \notag \\
       & = N \left |\sin \left (\frac \pi{N} \right) \right | \label{}
  \end{align}
  because  $\Delta x = 1/N$.
  %\geq (2 \pi S^*)^2 \sum_{k = -N/2, k \neq 0}^{N/2 -1} \left | \left ( \mathcal F_N \phi \right )_k \right|^2 (\Delta x)^2,
  Finally, using the Parseval identity \eqref{eq:DiscreteParseval} in \eqref{eq:almostPoincare} proves Lemma \ref{lem:poincare_discrete}, with $C_P = 1/S_*$.

\end{proof}

\begin{rem}
    In the continuous setting, the optimal constant for this Poincaré inequality is $C_P = 1/(2 \pi)$. In our case, when $N \to \infty$, $S_* \to \pi$, which is suboptimal. This is due to the fact that $D_x$ is a centered discretization of $\partial_x$. Taking a decentered discretization would yield $S_* \to 2 \pi$.
\end{rem}

\bibliographystyle{acm}
\bibliography{bib_BHR}

\end{document}